\documentclass[a4paper,leqno,11pt,final]{amsproc}
\usepackage{times,amssymb,latexsym,amsxtra,amscd,amsfonts,amsthm,amsmath,verbatim}

\begin{document}

\title[Signature Calculus]{Global Duality, Signature Calculus and the Discrete Logarithm Problem}

\author{Ming-Deh Huang }
\email{huang@pollux.usc.edu }
\address{Department of Computer Science\\ University of Southern California\\ Los Angeles\\ CA 90089-0781\\USA}

\author{Wayne Raskind}
\email{raskind@math.usc.edu }
\address{Department of Mathematics\\ University of Southern California\\ Los Angeles\\ CA 90089-2532\\ USA}

\newtheorem{theorem}{Theorem}
\newtheorem{lemma}{Lemma}
\newtheorem{corollary}{Corollary}
\newtheorem{definition}{Definition}
\newtheorem{remark}{Remark}
\newtheorem{conjecture}{Conjecture}
\newtheorem{proposition}{Proposition}
\newtheorem{algorithm}{Algorithm}
%below are from paper I
\newcommand{\bC}{{\Bbb C}}
\newcommand{\bZ}{{\Bbb Z}}
\newcommand{\bQl}{{\Bbb Q}_{\ell}}
\newcommand{\bZl}{{\Bbb Z}_{\ell}}
\newcommand{\Xb}{{\overline X}}
\newcommand{\Zb}{{\overline Z}}
\newcommand{\Eb}{{\overline E}}
\newcommand{\Kb}{{\overline K}}
\newcommand{\Cb}{{\overline C}}
\newcommand{\bGm}{{\mathbb G}_m}
\newcommand{\cX}{{\mathcal X}}
\newcommand{\Yb}{{\overline Y}}
\newcommand{\bQ}{{\mathbb Q}}
\newcommand{\bQp}{{\mathbb Q}_p}
\newcommand{\bZp}{{\mathbb Z}_p}
\newcommand{\bCp}{{\mathbb C}_p}
\newcommand{\bR}{{\mathbb R}}
\newcommand{\cE}{{\mathcal E}}
\newcommand{\bF}{{\mathbb F}}
\newcommand{\kb}{{\overline k}}
\newcommand{\Fb}{{\overline F}}
\newcommand{\cO}{{\mathcal O}}
\newcommand{\im}{\mbox{Im}}
\newcommand{\bZlZ}{{\bZ/\ell\bZ}}

\newcommand{\bFp}{{\mathbb F}_p}

\newcommand{\bFq}{{\mathbb F}_q}
\newcommand{\bFl}{{\mathbb F}_{\ell}}

\newcommand{\Z}{\mathbb Z}
\newcommand{\F}{\mathbb F}
\newcommand{\Q}{\mathbb  Q}
%above are from paper I

%below are from paper III
%\newcommand{\Z}{\mathbb Z}
%\newcommand{\F}{\mathbb F}
%\newcommand{\Q}{\mathbb  Q}
\newcommand{\E}{{\mathcal  E}}
\newcommand{\Fp}{\F_p}
\newcommand{\Fq}{\F_q}
\newcommand{\Fpn}{\F_{p^n}}

\newcommand{\co}{{\mathcal O}}

\newcommand{\coS}{{\mathcal O}_S}
\newcommand{\plim}{\displaystyle{\lim_{\stackrel{\longleftarrow}{n}}}\,}
\font\cyr=wncyr10
%above are from paper III

%\newcommand{\Kb}{{\overline K}}
\newcommand{\coK}{{\mathcal O}_K}
\newcommand{\coKu}{{\mathcal O}_{K_u}}
\newcommand{\coKv}{{\mathcal O}_{K_v}}

\maketitle

\begin{abstract} We study the discrete
logarithm problem for the multiplicative group and for elliptic
curves over a finite field by using a lifting of the corresponding
object to an algebraic number field and global duality. We introduce
the \textit{signature} of a Dirichlet character (in the
multiplicative group case) or principal homogeneous space (in the
elliptic curve case), which is a measure of the ramification at
certain places.  We then develop \textit{signature calculus}, which
generalizes and refines the index calculus method. Finally, we show
the random polynomial time equivalence for these two cases between
the problem of computing signatures and the discrete logarithm
problem.

\end{abstract}
\vspace{.5in}


\vspace{.5in}

\section{Introduction}
Let $A$ be a finite abelian group, which we write additively, and
let $x$ be an element of $A$. Let $y$ be in the subgroup generated
by $x$, so that $y=nx$ for some positive integer $n$.  Recall that
the {\it discrete logarithm problem} (DLP) is to determine $n$ in a
computationally efficient way. The computational complexity of this
problem when the bit size of the inputs is large is the basis of
many public-key encryption schemes used today. Two of the most
important examples of finite abelian groups that are used in
public-key cryptography are the multiplicative group of a finite
field and the group of points on an elliptic curve over a finite
field (see [Ko] and [Mill] for the original papers and [KMV] for a
survey of work as of 2000).\\

In what follows below, we will assume that $\ell$ is a large prime
number dividing the order of $A$ and that $x$ is an element of order
$\ell$. For $p$ a prime number and $q$ a power of $p$, we denote by
$\Fq$ the finite field with $q$ elements and by $\Fq^*$ its
multiplicative group of nonzero elements. \\

One of the best-known techniques to address the DLP is \textit{index
calculus}, which uses relations between elements of an abelian
algebraic group to derive linear relations between their discrete
logarithms. In the case of the multiplicative group of a finite
prime field, $\bFp$, taking sufficiently many random liftings of
elements of $\bFp^*$ to integers will ensure that  some  will only
be divisible by small (compared to $p$) prime numbers. Then such
relations can be derived because we know how to efficiently factor
integers that are products of powers of small prime numbers. See
e.g. [Mc], \S 5.1 or [SWD] for details. Trying to imitate this
method for an elliptic curve by lifting the curve to an algebraic
number field has turned out to be less effective, because the
behavior of the height function on the Mordell-Weil group of the
 lifted curve
 makes it far more difficult to derive relations like those just mentioned in the multiplicative group case  (see [HKT] or [JKSST] for more details).   However an important aspect of index calculus has not been
addressed in these studies, namely, the idea of leveraging small
primes to tackle a computational problem that involves large primes,
and it is not clear how this idea can be put to work in a setting
that involves the Mordell-Weil groups of elliptic curves. In this
paper we address this issue in both cases from the perspective of
arithmetic duality and propose a unified method which we call {\it
signature
calculus}.\\

Our general strategy to address the DLP in an abelian algebraic
group is to take a lifting of the group to an algebraic number field
and use the reciprocity law of global class field theory.  Others
have taken this approach (see e.g. [F] [FR], [N]), and we refine
their methods and give a general exposition of the theory.  We
explain below in detail how this works for the multiplicative group
of a finite field and for the group of points of an elliptic curve
over a finite field. The idea is to construct a suitable ``test''
element, which is a Dirichlet character in the multiplicative group
case and a principal homogeneous space in the elliptic curve case.
This element pairs with the lifting of a point of the group to give
an equation between the local terms of this pairing. The lifting
from a finite field $\F_p$ to a global field preserves discrete
logarithms at a place over $p$.  The reciprocity law then allows us
to distribute information on the discrete logarithms among a set of
places which depends on the choice of test element and the manner of
lifting. We define the {\it signature} of these test elements and
prove the equivalence of computing the signature with the respective
DLP. These signatures measure the ramification at primes above $p$
and $\ell$.   Though the signatures are small, they uniquely
identify the objects they represent (Dirichlet characters and
principal homogeneous spaces). They are, in fact, succinct
representations of those objects, and the equivalence results show
that computing these signatures (without constructing the objects
they succinctly represent) amounts to solving discrete-log
problems.\\

The unifying approach based on global duality provides an ideal
setting to compare and contrast index calculus methods in the
multiplicative group and elliptic curve cases. The signature
computation problem involves large primes, and the question
naturally arises as to whether small primes can be utilized to
tackle the problem with greater computational efficiency, in a
similar way as we mentioned for the multiplicative group. Following
the equivalence results we show that in this setting, the index
calculus method arises quite naturally for the discrete-log problem
in the multiplicative case and the corresponding signature
computation problem.  In contrast, a similar method cannot be
fashioned for the elliptic curve case. The success in one case and
the lack thereof in the other is due to the difference in the nature
of the pairings involved.  In the multiplicative case, a Dirichlet
character which is unramified at a finite place $v$ can nevertheless
pair nontrivially with local non-units at $v$. This makes it
possible for small primes to play a role in forming relations among
values of local pairings. In the elliptic curve case, an unramified
principal homogenous space at a good reduction place $v$ is one that
extends to a principal homogeneous space under a smooth proper model
$\mathcal{E}_v$ of $E$ over the ring of local integers $R_v$ (please
see \S 1 below for more details and explanation). There is a
bijection between such principal homogeneous spaces and the
corresponding objects under the reduction of $\mathcal{E}_v$ mod $v$
(see e.g. [MET], Chapter III, Remark 3.11(a)).  By a theorem of Lang
([L], Theorem 2), the latter objects are trivial.  Thus, in the
elliptic curve case, an unramified principal homogeneous space at a
good reduction place is trivial. For small primes of bad reduction
not dividing $\ell$, only the group of components of the special
fibre of the N\'eron model of the elliptic curve over the ring of
integers plays a role, and the order of this group is unlikely to be
divisible by $\ell$ (see \S 5.1.2 below for more details). As a
result, only primes of large norm can play a role in forming
relations among values of local pairings in the elliptic
curve case.\\

The computational complexity of signature calculus is an intriguing
question, since the objects involved (Dirichlet characters and
principal homogeneous spaces) and their associated field extensions
are huge, but the signatures sought are small. Although we show that
the testing Dirichlet characters and principal homogeneous spaces
exist, it remains an interesting question as to how they can be
explicitly constructed. This is easier to handle in the
multiplicative case, where we also derive a concrete number
theoretic characterization of the character signature by working out
the local pairings using norm residue symbols. For the elliptic
curve case, we have a partial
solution. \\

This paper is a more formal and detailed exposition than the survey
of this material that appeared in [HRANTS], and it contains very
significant material that is not in that paper.  We have tried to be
completely mathematically precise while retaining the cryptographic
motivation and applications. \\

The idea of using global methods in this way was originally proposed
by Frey [F], whom we thank for inspiration, helpful discussions, and
for inviting us to present our work at the Elliptic Curve
Cryptography (ECC) conference in Bochum in September 2004. Methods
of this type have also been used by Frey and R\"uck [FR], and by
Nguyen [N].
\section{Global Framework}
\subsection{Notation and Preliminaries}
If $A$ is a locally compact abelian group that is either profinite
or torsion, we denote by $A^*$ the group
$\mbox{Hom}_{cont}(A,\Q/\Z)$ of continuous homomorphisms and refer
to it as the {\it Pontryagin dual} of $A$.  Note that * is an
exact functor since $\Q/\Z$ is a divisible abelian group.\\

Let $K$ be a field, fix a separable closure $\Kb$ of $K$, and let
$G=Gal(\Kb/K)$.  Let $M$ be a discrete $G$-module upon which $G$
acts continuously, where $G$ has the Krull topology.  We will be
using Galois cohomology extensively, which we will denote by
$H^i(G,M)$ or sometimes $H^i(K,M)$ . A basic reference for this
theory is [S1].\\

We shall mainly be using three types of fields: finite fields,
denoted by $\mathbb F$, algebraic number fields, denoted by $K$, and
the
completion of an algebraic number field at a finite place $v$,  denoted by $K_v$.\\

 An \textit{algebraic number field} will be a finite
extension of the field of rational numbers $\bQ$. We consider
equivalence classes of absolute values $v$ on $K$, which we call
{\it places}. As most of our discussion will pertain to abelian
groups that are $\ell$-torsion, where $\ell$ is an odd prime number,
we shall
ignore the real and complex places for the most part.\\

Let $R$ be a discrete valuation ring with fraction field $K$ and
residue field $F$.  For example, $R$ could be the ring of integers
in a $K_v$. Let $X$ be a smooth proper scheme over
$Y=\mbox{Spec}(R)$. Recall that this means that the structure
morphism:

$$f:X\to Y $$
is smooth and proper.  The former condition means that the fibres
over $K$ (the generic fibre) and $F$ (the special fibre) are smooth,
and the latter means that $f$ is separated and universally closed
(i.e. that if we change base by a morphism $Z\to Y$, then the
morphism:

$$X\times_YZ\to Z$$
is closed).  If $X\to Y$ is a proper morphism, then a point $P\in
X(K)$ may be lifted to a point in $X(R)$.  If $E$ is an
 elliptic curve over $K$, we may clear the denominators in a
 defining equation and view it as a curve over $R$ (not necessarily smooth over $R$).  Then
 $E$ is proper over $R$, whereas the multiplicative group is affine
 and decidedly not proper.\\

Recall that an \textit{elliptic curve} over a field $K$ is a smooth,
projective algebraic curve $E$ of genus 1 together with a
distinguished rational point {\it O}, which serves as the identity
element in an abelian group structure on $E$ that can be defined
geometrically by a chord and tangent method.  We denote  by $E(K)$
the set of points of $E$ over $K$. Recall that a {\it principal
homogeneous space} under $E$ over $K$ is a curve $F$ of genus 1 over
$K$ together with a simply transitive group action of $E$ on $F$.
The isomorphism classes of such principal homogeneous spaces are
classified by the group $H^1(G,E(\Kb))$, where $G=Gal(\Kb/K)$. A
principal homogeneous space is trivial if and only if it has a
rational point over $K$, in which case it is isomorphic to $E$ over
$K$. Thus any principal homogeneous space becomes isomorphic to $E$
over a finite extension of $K$.\\

Let $\mathcal{M}$ be an algebraic group over a discrete valuation
ring $R$ and denote by $M$ its fibre over the quotient field $K$. We
will be most interested in the cases where $\mathcal{M}$ is either
the constant algebraic group $\Z/\ell\Z$ or a smooth proper model of
an elliptic curve $E$ with good reduction over a completion of an
algebraic number field at a finite place, $v$. Recall that an
element of $H^1(K,M)$ is said to be \textit{unramified} if it is in
the image of the natural map:

$$H^1(R,\mathcal{M})\to H^1(G,M).$$

This is a more general notion of non-ramification, which is the same
as the usual definition when $M$ is finite.\\

 Let $\tilde{E}$ be an elliptic curve over $\F$ and let $R$ be a discrete valuation ring
$R$ with quotient field $K$ and residue field $\bF$. Then a {\it
lifting} $E$ of $\tilde{E}$ to $K$ is a smooth proper scheme $\cE$
over $R$ whose special fibre is $\tilde{E}$ and whose generic fibre
is $E$. We shall use rather simple liftings below, but let us point
out that it is a theorem of Deuring [D] that if $\tilde{E}$ is an
elliptic curve
 over a finite field  with an endomorphism $\varphi$, then the pair $(\tilde{E},\varphi)$ can be lifted to
 to a discrete
valuation ring $R$ whose quotient field is an algebraic number
field.  If the curve is \textit{ordinary}, as are the curves we
consider in this paper, then one can lift the curve together with
the whole endomorphism ring.  A more systematic approach to
liftings of ordinary elliptic curves is given by Serre-Tate theory (see e.g. [S2], \S 5).\\

Recall the Brauer group $Br(K)$ of similarity classes of finite
dimensional central simple algebras over $K$, which can be described
 in terms of Galois cohomology by

$$Br(K)\cong H^2(G,\Kb^*).$$

When $K$ is an algebraic number field, we have the
Brauer-Hasse-Noether exact sequence:

$$(\dag)\:0\to Br(K)\to \sum_vBr(K_v)\to \bQ/\bZ\to 0.$$
This is the beginning of the theory of {\it global duality}, which
shows how to relate the arithmetic of $K$ with that of all of the
$K_v$. The following subsections review this theory briefly in the
context in which we shall use it.

\subsection{Reciprocity Law for the Multiplicative Group}  We review the reciprocity law in this
context, mostly following the exposition of ([S1], Chapter XIV). Let
$K^*$ denote the set of nonzero elements of $K$, which is an abelian
group under multiplication. We consider a Dirichlet character $\chi$
of $K$, which we view as an element of the Galois cohomology group
$H^1(G,\bQ/\bZ)$.  Thus $\chi$ represents a finite cyclic extension
$L/K$ together with a homomorphism:

$$Gal(L/K)\to  \bQ/\bZ.$$

Let $\partial(\chi)$ denote the image of $\chi$ under the boundary
map

$$H^1(G,\bQ/\bZ)\stackrel{\partial}{\to} H^2(G,\bZ)$$
in the long exact cohomology sequence associated to the short exact
sequence of $G$-modules with trivial action:

$$0\to \bZ\to \bQ\to\bQ/\bZ\to 0.$$
Then for $a\in K^*$ we consider

$$<\chi,a>:=a\cup\partial(\chi) \in H^2(G,\Kb^*)$$
under the pairing:
$$K^*=H^0(G,\Kb^*)\times H^2(G,\bZ)\to H^2(G,\Kb^*)\cong Br(K).$$

If $L$ is the extension corresponding to $\chi$, then we have that
$<\chi,a>=0$ if and only if $a$ is a norm from $L^*$.\\

    If $K$ is an algebraic
number field, $\chi\in H^1(G,\bQ/\bZ),  a\in K^*$ and $v$ is a
 place of $K$, then we can restrict $\chi$ to each
$K_v$ and regard $a$ as an element of $K_v^*$.  Note that we may
have $\chi_v=0$. We then denote the local pairing by $<\chi_v,a_v>$.
If $v$ is a nonarchimedean place then $Br(K_v)\cong \bQ/\bZ$ and we
view $<\chi_v,a_v>$ as an element of $\bQ/\bZ$.  Note also that if
$v$ is a place where $\chi$ is unramified and $a$ is a unit at $v$,
then $<\chi_v,a_v>=0$.  That is, every unit is a norm from an
unramified extension of nonarchimedean local fields.   Thus
$<\chi_v,a_v>$=0 for all but finitely many $v$. Since the local
pairings are compatible with the global pairings, the exact sequence
($\dag$) above for the Brauer group of an algebraic number field
shows that we have the {\it reciprocity law}

$$\sum_v <\chi_v,a_v>=0\in \bQ/\bZ.$$

\subsection{Reciprocity Law for Elliptic Curves}
Let $E$ be an elliptic curve over $K$.   Let $Q\in E(K)$ and
$\alpha\in H^1(K,E)$. We consider the pairings

$$<\alpha,Q>\in Br(K)$$
$$<\alpha_v,Q_v>\in Br(K_v)\cong\bQ/\bZ .$$

These are not as easy to describe explicitly as in the case of the
multiplicative group, but we give here a quick if somewhat terse
definition.  Given an abelian variety $A$ over $K$, let $\hat{A}$
denote its dual, which is $Ext_K^1(A,\bGm)$, where $\bGm$ is the
multiplicative group scheme and the $Ext$ is taken in the category
of algebraic groups over $K$. An elliptic curve is self-dual, so
that we can identify $E(K)$ with $Ext_K^1(E,\bGm)$.  Given $Q\in
E(K)$, represent it as a 1-extension of algebraic groups using this
identification

$$0\to \bGm\to X\to E\to 0,$$
and let

$$(\dag\dag)\: 0\to \Kb^*\to X(\Kb)\to E(\Kb)\to 0$$
be the short exact sequence of $\Kb$-points of these groups. Then
given an element $\alpha\in H^1(G,E(\Kb))$, let
$<\alpha,Q>=\partial_Q(\alpha)$, the image of $\alpha$ under the
boundary map:

$$H^1(G,E(\Kb))\stackrel{\partial_Q}{\to}H^2(G,\Kb^*)$$
in the long exact cohomology sequence obtained from the short exact
sequence ($\dag\dag$). For $\alpha\in H^1(G,E(\Kb))$ and $Q\in E(K)$
we denote by $\alpha_v$ the image of $\alpha$ in $H^1(G_v,E(\Kb_v))$
(which may be zero) and by $Q_v$ the image of $Q$ in $E(K_v)$.
 We
can make a similar definition over the nonarchimedean fields $K_v$
for $\alpha_v\in H^1(G_v,E(\Kb_v))$ and $Q_v\in E(K_v)$ to get
$<\alpha_v,Q_v>\in Br(K_v)\cong \bQ/\bZ$. \\

We will be interested in the situation where $\alpha\in H^1
(K,E)[\ell]$, in which case we have the following commutative
diagram:

\[
\begin{array}{clclc}
E(K)/\ell & \times & H^1 (K, E )[\ell]   & \rightarrow & Br (K)[\ell] \\
\downarrow &  & \downarrow & & \downarrow\\
E(K_v )/\ell & \times & H^1 (K_v , E)[\ell] & \rightarrow & Br (K_v)
[\ell]
\end{array}
\]

In the case of the local field $K_v$, the pairing is perfect (local
duality for abelian varieties, see e.g. [MAD], Ch. I, \S 3,
Corollary 3.4). \\

We then have that $<\alpha_v,Q_v>=0$ for almost all $v$. The
fundamental sequence ($\dag$),
%$$0\rightarrow Br(K) \rightarrow \oplus_v  Br(K_v)
%\stackrel{\Sigma_v {\rm inv}_v}{\rightarrow} {\Q}/{\Z}\rightarrow
%0,$$
the identification $Br(K_v)[\ell]\cong \Z/\ell\Z$, and the
commutative diagram above imply that for $\alpha\in H^1 (K,E
)[\ell]$ and $Q\in E(K)$,

$$\sum_v <\alpha_v,Q_v>=0\in \bQ/\bZ.$$

\subsection{Cohomological Basis of the Unified Approach}
Our approach is based on duality theorems for Galois modules and for
abelian varieties over number fields. Let $K$ be an algebraic number
field and $\cO_K$ the ring of integers in $K$.
 Let $X=\mbox{Spec}(\cO_K)$ and $U$ be a nonempty open subset
of $X$ with complement $S$. Thus $U$ consists of all but finitely
many places of $K$.  Let $\ell$ be a prime number that is invertible
on $U$ and let $\mu_{\ell}$ be the sheaf of $\ell$-th roots of
unity.  We are interested in the groups $H^i(U,\mu_{\ell})$. To aid
us in computing them and related cohomology groups, we have the {\it
Poitou-Tate exact sequence} (see e.g. [MAD], Ch. I, \S 4, Theorem
4.10c):

$$0\to H^0(U,\mu_{\ell})\to \bigoplus_{v\in S}H^0(K_v,\mu_{\ell})\to
H^2(U,\bZ/\ell\bZ)^*\to$$
$$H^1(U,\mu_{\ell})\to \bigoplus_{v\in
S}H^1(K_v,\mu_{\ell})\to H^1(U,\bZ/\ell\bZ)^*\to$$
$$H^2(U,\mu_{\ell})\to \bigoplus_{v\in S}H^2(K_v,\mu_{\ell})\to
H^0(U,\bZ/\ell\bZ)^*\to 0.$$

This sequence summarizes many of the basic results from class field
theory. Let $K_S$ be a maximal extension of $K$ that is unramified
outside $S$ and put $G_S=Gal(K_S/K)$. Then any sheaf $\mathcal{F}$
on $U$ may be regarded as a $G_S$-module, and we have
$H^i(U,\mathcal{F})\cong H^i(G_S,\mathcal{F})$. We shall often use
this latter notation for the multiplicative group case.  We are
mainly interested in the middle line of the Poitou-Tate sequence:

$$(*)_{\mu_{\ell}}:\:H^1(G_S,\mu_{\ell})\to \bigoplus_{v\in
S}H^1(K_v,\mu_{\ell})\to H^1(G_S,\bZ/\ell\bZ)^*$$ and the dual
sequence obtained by taking the Pontryagin dual and using Tate local
duality:

$$(*)_{\Z/\ell\Z}\::H^1(G_S,\bZ/\ell\bZ)\to \bigoplus_{v\in
Z}H^1(K_v,\bZ/\ell\bZ)\to H^1(G_S,\mu_{\ell})^*.$$

 For an elliptic curve $E$ over $K$ that has a
smooth proper model ${\mathcal E}$ over $U$ on which $\ell$ is
invertible, we have the {\it Cassels-Tate exact sequence} (see
[MAD], Ch. II, \S 5, Theorem 5.6b):

$$(**)\:  E(K)^{(\ell)}\to \bigoplus_{v\in S}E(K_v)^{(\ell)}\to H^1(U,{\mathcal
E})\{\ell\}^*\to \hbox{\cyr Sh}(E)\{\ell\}\to 0.$$

Here $(\ell)$ denotes completion with respect to subgroups of
$\ell$-power index, $\{\ell\}$ denotes the $\ell$-primary part of a
torsion abelian group, and $\hbox{\cyr Sh}(E)$ is the
Shafarevich-Tate group of
everywhere locally trivial principal homogeneous spaces under $E$, which we assume to be finite.\\

We give here a very terse explanation of the common origin of these
two exact sequences, as it is the key to our unified approach in the
multiplicative group and elliptic curve cases. Let $\mathcal{F}$ be
a sheaf on $U$ and $j_!\mathcal{F}$ denote extension of
$\mathcal{F}$ by zero from U to $X$.  We denote by
$H^i_c(U,\mathcal{F})$ the group $H^i(X,j_!\mathcal{F})$; this is
{\it cohomology with compact support}. Then we have a long exact
sequence of cohomology with support (see [MET], Chapter III,
Proposition 1.25):

$$\cdots H^i_S(X,j_!\mathcal{F})\to H^i(X,j_!\mathcal{F})\to H^i(U,j^*j_!\mathcal{F})\to
H^{i+1}_S(X,j_!\mathcal{F}).$$

For a place $v$ of $K$, let $A_v^h$ denote the henselization of the
local ring of $X$ at $v$ (one can also take the completion).  Then
using the identifications:

$$H^i_S(X,j_!\mathcal{F})\cong \bigoplus_{v\in S}H^i_v(X,j_!\mathcal{F})$$

$$H^i_v(X,j_!\mathcal{F})=H^i_v(A_v^h,j_!\mathcal{F})$$

$$H^i(K_v,\mathcal{F})\cong H^{i+1}_v(A_v^h,j_!\mathcal{F})$$
%\vspace{.1in}
for $v\in S$ (see [MAD], Proposition 1.1, page 182 for the last
isomorphism, which uses the fact that we have a sheaf of the form
$j_!\mathcal{F}$), we get the exact sequence

$$\cdots H^i_c(U,\mathcal{F})\to H^i(U,\mathcal{F})\to \bigoplus_{v\in
S}H^i(K_v,\mathcal{F})\to H^{i+1}_c(U,\mathcal{F})\cdots.$$

 The
Poitou-Tate and Cassels-Tate exact sequences are then derived from
this one sequence by taking $\mathcal{F}=\mu_{\ell}$ (resp.
$\mathcal{F}={\mathcal E}$) and using the Artin-Verdier duality
theorem (see e.g. [MAD], Chapter II, \S 3, Corollary 3.2) (resp. the
duality theorem for abelian
varieties (see [MAD], Chapter 3, \S 5, Theorem 5.2)).\\

\section{Classical Index Calculus from the Perspective of
Arithmetic Duality}

Our approach to the discrete log problem for the multiplicative
group of a finite field uses the Poitou-Tate exact sequence (*) in
\S 2 above. For the discrete log problem for an elliptic curve
$\tilde{E}$ over a finite field with a point of order $\ell$ and a
suitable lifting $E$ of $\tilde{E}$ to an algebraic number field
$K$, we will use the Cassels-Tate sequence (**) in \S 2, where $U$
is an open subset of $\mbox{Spec}(\cO_K)$ on which $E$ has good
reduction and $\ell$ is invertible, and ${\mathcal E}$ is a smooth
proper model of $E$ over $U$.  In each case, the method will be to
find a suitable element of $H^1(U,\mathcal{F})$ of order $\ell$
against which to ``test'' a lifting to $K$ of an element over the
finite field whose discrete log we seek to compute and then use the
reciprocity laws that are encoded in the exact sequences to create
linear
relations between the discrete logs.\\

We demonstrate below how the classical index calculus method emerges
in this context as the result of one
particular choice of testing Dirichlet character and method of lifting.\\

Let $p$ and $\ell$ be odd primes such that $p\equiv 1 \pmod{\ell}$
but $p\not\equiv 1 \pmod{\ell^2}$. Given positive integers $g$ and
$t$ such that $g\bmod{p}$ generates the group $\F_p^*$, we would
like to compute $n\bmod{\ell}$ where $t = g^n$ in $\F_p^*$.  We will
fix $g$ and denote the discrete-log $t$ with respect to $g$ as
$\theta (t)$. The core of the classical index calculus method for
solving the discrete-log problem in $\F_p^*$ is to compute
$\theta(q)$ for primes $q$ up to a chosen bound $B$.\\

Let $K=\Q$, $X=\mbox{Spec}(\Z)$, and $U=X-S$, where $S$ is a finite
set of primes containing $\ell$. Consider the sequence
$(*)_{\Z/\ell\Z}$ of the last section.
 The extension $\Q(\mu_p)/\Q$ is cyclic of degree
$p-1$. Since $p\equiv 1 \pmod{\ell}$, there is a unique
sub-extension $L/\Q$ of degree $\ell$.  We fix an isomorphism
$Gal(L/\Q)\cong \Z/\ell\Z$ and denote by $\chi$ the corresponding
Dirichlet character, which is ramified only at $p$. Then $\chi$ can
be regarded as an element of $H^1(G_S,\Z/{\ell}\Z)$ if $p\in S$. We
have that $\Z_S^*/{\Z}_S^{*\ell}\cong H^1(G_S,\mu_{\ell})$, and from
$(*)_{\Z/\ell\Z}$ we have that for all $\alpha\in\Z_S^*$,
$$\sum_{v\in S} < \chi_v,\alpha_v >=0\in \Z/\ell\Z.$$
Note that
$$ < \chi_p,\alpha_p > = \theta(\alpha) <\chi_p,g > ,$$
and for $q\in S-\{p\}$,
$$ < \chi_q,\alpha_q > = v_q (\alpha) <
\chi_q,q >,$$ where $v_q (\alpha)$ is the $q$-adic valuation of
$\alpha$.\\

Let $F$ be the set of primes up to some bound $B$ and let $S$ be the
set $F$ together with $p$ and $\ell$. For $q\in F$, since
$q\in\Z_S^*$ and $q$ is a local unit at $v\neq q$ in $S$,
$$0= \sum_{v\in S} < \chi_v,q > = < \chi_p ,q > + < \chi_q ,q >
= \theta(q) < \chi_p, g > + < \chi_{q}, q > .$$ Hence,
$$\theta(q) = -(< \chi_p, g >)^{-1}  < \chi_{q}, q >.$$

To compute $\theta(q)$ for all primes $q$ in $F$, we generate random
$r$ so that $g^r\bmod{p}$ is $B$-smooth, that is
$$\alpha_r = g^r \bmod{p} = \prod_{q\in F} q^{e_q (r)}$$
with $e_q(r)\in\Z_{\ge 0}$.  Since $\alpha_r \in \Z_S^*$, we have
$$0= \sum_{v\in S} < \chi_v , (\alpha_r)_v > = r < \chi_p, g > +
\sum_{q\in F} e_q(r) < \chi_q, q >.$$ It follows that
$$\sum_{q\in F} e_q (r) \theta(q) = r.$$
With sufficiently many $\alpha_r$ that generate
$\Z_F^*/\Z_F^{*\ell}$, we can solve for the unknown
$\theta(q)\bmod{\ell}$.
What we have derived is in essence the classical index calculus method.\\

We remark that similar reasoning as above shows that the image of
$H^1 (G_S, \Z/\ell\Z)$ in $\bigoplus_{v\in S}H^1(\Q_v,\Z/{\ell}\Z)$,
where $S=F\cup\{p\}$, has $\F_{\ell}$ dimension one, and the
classical index calculus method amounts to
determining this image in a computationally efficient manner.\\

In the preceding discussion, we were able to explicitly construct
the desired Dirichlet character because we were working with abelian
extensions of  $\Q$, about which we know enough to explicitly
compute everything we need. In the discussion below we will be
working with real quadratic fields, and there we know much less
about how to explicitly construct abelian extensions. However, using
the exact sequence $(*)_{\mu_{\ell}}$, we will demonstrate the
existence of a suitable Dirichlet character by explicitly computing
the $\bFl$-dimensions of the first and second terms, and showing
that the former is less than the latter.  More generally, we use the
following basic strategy to find a suitable testing element. In the
multiplicative group case, look for an algebraic number field $K$
such that the $\bFl$-dimension of the first term of the middle row
of $(*)_{\mu_{\ell}}$ is smaller than that of the second. This will
then guarantee the existence of an element of order $\ell$ in
$H^1(G_S,\Z/\ell\Z)^*$. By lifting to units of a real quadratic
field instead of to smooth integers in $\Z$, we are more able to
compare and contrast the discrete log problems for the
multiplicative group and for elliptic curves over finite fields.  It
is an artifact of class field theory that one can often demonstrate
the existence of an abelian extension without there being an obvious
way to construct it explicitly.\\

In the elliptic curve case, we look for an algebraic number field
$K$ together with an elliptic curve $E/K$ that lifts $\tilde{E}$,
such that $E(K)$ is of small rank, e.g. $\leq 2$.  We also assume
that at least one of the generators of the torsion-free quotient of
$E(K)$ is not divisible by $\ell$ in $E(K_u)$ for all $u\in T$,
where $T$ consists of one place above $p$ and both above $\ell$ in a
quadratic extension $K/\bQ$ in which both
$p$ and $\ell$ split.\\

This approach will be developed in more detail in the next few
sections.

\section{Signature Calculus for the Multiplicative Group}
\subsection{Characters with Prescribed Ramification}
\label{character}
 Throughout this section, let $p,{\ell}$ be
rational primes with $p\equiv 1 \pmod{{\ell}}$ and $\ell > 2$. Let
$K/\Q$ be a real quadratic extension where $p$ and ${\ell}$ split.
Let $\alpha$ be a fundamental unit of $K$. Let $\Sigma$ be the set
of all places over $\ell$ and $p$, together with all the archimedean
places. For any place $u$ of $K$ let $P_u$ denote the prime ideal
corresponding to $u$. For any finite set $S$ of places of $K$, let
$G_S$ denote the Galois group of a maximal extension of $K$ that is
unramified outside of $S$.

\begin{proposition}
\label{ram} Let $S$ be a subset of $\Sigma$ that contains both
places over $\ell$ and both archimedean places.  Suppose
\begin{enumerate}
\item $\ell \nmid h_K$ where $h_K$ is the class number of $K$;
\item either $\alpha^{l-1} \not\equiv 1 \pmod{P_w^2}$  for some
$w\in S$ over $\ell$, or $\alpha^{\frac{p-1}{{\ell}}} \not\equiv 1
\pmod{P_w}$ for some $w\in S$ over $p$ (that is, locally $\alpha$ is
not an $\ell$-th power at either a place over $\ell$ or a place over
$p$).
\end{enumerate}
Then the $\F_{\ell}$-dimension of  $H^1 (G_S , \Z/{\ell}\Z )$ equals
$n(S)-1$ where $n(S)$ is the number of finite places in $S$.
\end{proposition}

\ \\{\bf Proof}:   Consider the sequence:

$$(*)_{\mu_{\ell}}:\:H^1(G_S,\mu_{\ell})\stackrel{f}{\to} \bigoplus_{v\in
S}H^1(K_v,\mu_{\ell})\stackrel{\rho}{\to} H^1(G_S,\bZ/\ell\bZ)^*\to
H^2(G_S,\mu_{\ell})\stackrel{g}{\to} \bigoplus_{v\in
S}H^2(K_v,\mu_{\ell})\to\cdots.$$

We claim that under the hypotheses of the proposition, $\rho$ is
surjective.  To see this, the hypothesis that $\ell$ does not divide
the class number of $K$ implies that it does not divide the class
number of ${\mathcal O}_S$.  By Kummer theory, we then have that:

$$H^2(G_S,\mu_{\ell})\cong Br({\mathcal O}_S)[\ell].$$

But then the map $g$ is injective, so $\rho$ is surjective.   Now
consider the map

$$f:\:H^1(G_S,\mu_{\ell}) \stackrel{f}{\to} \bigoplus_{v\in S}H^1(K_v,\mu_{\ell}).$$

Again using the hypothesis that $\ell$ does not divide the class
number of $K$, we have that:

$${\mathcal O}_S^*/{\mathcal O}_S^{*\ell}\cong H^1(G_S,\mu_{\ell}).$$

Consider the exact sequence:

$$0\to {\mathcal O}^*\to {\mathcal O}_S^*\to \Z S\to Cl({\mathcal
O}) \to Cl({\mathcal O}_S)\to 0.$$

Going modulo $\ell$ and using the hypotheses of the theorem, we see
that the sequence:

$$0\to {\mathcal O}^*/{\mathcal O}^{*\ell}\to {\mathcal O}_S^*/{\mathcal O}_S^{*\ell}\to \Z S/\ell\Z S\to
0$$
is exact.  This shows that the $\F_{\ell}$-dimension of the
group in the middle is $n(S)+1$.  The hypotheses about the units
show that $f$ is injective. The target has dimension $2n(S)$ because
$H^1(K_v,\mu_{\ell})$ is isomorphic to $\Q_v^*/\Q_v^{*\ell}$.  If
$v\mid p$, then this group is of dimension 2 over $\F_{\ell}$
because $\ell\mid p-1$.  If $v\mid \ell$, then this group is also of
dimension 2, spanned by a prime element of $\Q_{\ell}$ and by a
1-unit.  Thus the cokernel of $h$ is of dimension $n(S)-1$.  This
completes the proof of the proposition.

\begin{proposition}
\label{oneram} Let $S$ be the set consisting of one place $u$ over
$\ell$, one place $v$ over $p$, and both archimedean places. Suppose
\begin{enumerate}
\item $\ell \nmid h_K$ where $h_K$ is the class number of $K$;
\item $\alpha^{l-1} \not\equiv 1 \pmod{P_w^2}$ for all places $w
\mid \ell$; \item $\alpha^{\frac{p-1}{{\ell}}} \not\equiv 1
\pmod{P_v}$.
\end{enumerate}
Then the $\F_{\ell}$-dimension of $H^1 (G_S , \Z/\ell \Z )$ is one.
If $\chi$ is any nonzero element of this group, then $\chi$ is
ramified at $u$ and $v$.
\end{proposition}
\ \\{\bf Proof} Suppose $u,u'$ are the places over $\ell$. Let $R$
be the set consisting of $u,u'$ and both archimedean places.  Let
$T$ be the set consisting of $u,u', v$ and both archimedean places.
Then from Proposition~\ref{ram} it follows that $H^1 (G_R ,\Z/\ell
\Z )$ has dimension one and $H^1 (G_T ,\Z/\ell \Z )$ has dimension
two.  Hence there exists a nontrivial $\psi\in H^1 (G_R ,\Z/\ell \Z
)$, and some $\chi\in H^1 (G_T ,\Z/\ell \Z ) - H^1 (G_R ,\Z/\ell \Z
)$. By construction $\chi$ is ramified at $v$, and by the condition
on $\alpha$ at $v$ we get $< \chi_v , \alpha_v
>\neq 0$. As for $\psi$, by the reciprocity law we have
$< \psi_u , \alpha_u > + < \psi_{u'} , \alpha_{u'} > =0$, so either
$< \psi_u , \alpha_u >$ and  $< \psi_{u'} , \alpha_{u'} >$ are both
zero or both non-zero.  But if both are zero then by the condition
of $\alpha$ at $u$ and $u'$ it would follow that $\psi$ is
unramified at both places, violating the condition that $\ell$ does
not divide the class number of $K$. Hence $< \psi_u , \alpha_u >$
and  $< \psi_{u'} , \alpha_{u'} >$ are both non-zero.  Since $<
\psi_{u'} , \alpha_{u'}
> \neq 0$, there exists $c\in\Z/\ell\Z$ such that $ < \chi_u ,
\alpha_{u'}
> = c < \psi_{u'} ,\alpha_{u'} > $, and letting $\phi = \chi - c\psi$, we have
$ < \phi_{u'} , \alpha_{u'}
> = 0$.
Now $\phi\in H^1 (G_S,\Z/\ell\Z )$ since  $ < \phi_{u'} ,
\alpha_{u'}
> = 0$, and $\phi$ is a nontrivial since $< \phi_v , \alpha_v >
= < \chi_v , \alpha_v > \neq 0$.  Hence $H^1 (G_S,\Z/\ell\Z )$ is of
dimension at least one.  Since $\psi$ is ramified at $u'$, it
follows that $\psi\not\in H^1 (G_S ,\Z/\ell \Z )$, and since
$\psi\in H^1 (G_R ,\Z/\ell \Z )\subset H^1 (G_T ,\Z/\ell \Z )$, it
follows that $H^1 (G_S ,\Z/\ell \Z )$ is a proper subset of $H^1
(G_T ,\Z/\ell \Z )$, hence it can be of dimension at most one.  We
conclude that its dimension must be one, and the proposition
follows.\\

\subsection*{Remarks}

\noindent (i)  We explain why we made the assumptions of Proposition
2, their necessity and sufficiency
for the conclusion, and how they affect the signature computations later in the paper:\\

Condition (1) is made to ensure that the Dirichlet characters of
degree $\ell$ that we get will not be everywhere unramified, as
such characters would be of no use to us for the signature computation.\\

Conditions (2) and (3) are meant to ensure that there do not exist
characters of $K$ of degree $\ell$ that are ramified only at
$\mathfrak p$ or only at $\mathfrak l$.  Such characters would not
help our signature computation. For example, suppose the character
$\chi$ is ramified at $v$ and unramified everywhere else. Then if we
pair $\chi$ with a global unit $a$ of our real quadratic field, we
would get that $<\chi_u,a_u>=0$ since $\chi$ is unramified at $u$
and $a$ is a unit.  The reciprocity law would then give us that
$<\chi_v,a_v>=-<\chi_u,a_u>=0$, and this would not help us in the
signature computation.  If {\it neither} condition (2) nor (3)
holds, then there are Dirichlet characters $\chi'$ and $\chi''$, one
ramified only at $u$ and the other ramified only at $v$. Thus, while
the character $\chi=\chi'+\chi''$ is ramified at both $u$ and $v$,
this would not help for our signature computation, since for a
global unit $a$, we would have:

$$<\chi_u,a_u>=<\chi'_u,a_u>+<\chi''_u,a_u>=0,$$
since $\chi'$ is ramified only at $u$ and $\chi''$ is unramified at
$u$.  Similarly for $v$.\\

(ii)  One can give an alternative (and perhaps simpler) proof of
Proposition 2 using the ideal theoretic formulation of class field
theory. Briefly, using the hypotheses of the proposition, one easily
calculates the $\ell$-rank of the Galois group of the ray class
field modulo $I=\mathfrak p\mathfrak l^2$, where $\mathfrak p$ is an
ideal of $K$ lying over $p$ and $\mathfrak l$ is an ideal lying over
$\ell$. This is the maximal abelian extension of $K$ with conductor
bounded by $I$, and its Galois group is isomorphic to a generalized
class group by class field theory.  Using basic exact sequences and
the hypotheses of the proposition, we can explicitly calculate this
class group.  The reason why we did not write the proof this way is
that we want to stress the analogy with elliptic curves, where the
Poitou-Tate exact sequence has an analogue (the Cassels-Tate exact
sequence), but the ideal theoretic formulation of
class field theory has no known analogue.\\

Assuming the conditions in Proposition~\ref{oneram}, then $H^1
(G_S,\Z/{\ell}\Z)$ is isomorphic to $\Z/\ell\Z$.  Every nontrivial
character in it is ramified at $u$ and $v$ and unramified at all
other finite places; moreover, $<\chi_u , \alpha_u
> \neq 0$ and $<\chi_v , \alpha_v >\neq 0$, and $<\chi_u ,
\alpha_u
> + <\chi_v , \alpha_v > = 0$.
This group of characters corresponds to a unique cyclic extension
$K_S$ of degree ${\ell}$ over $K$ which is ramified at $u$ and $v$
and unramified at all other finite places.\\

At $u$, we take the class of $1+\ell$ as the generator of the group
$\coKu^* /\coKu^{*\ell} \cong \Z_{\ell}^* /\Z_{\ell}^{*\ell} \cong
\Z/{\ell}\Z$. For $0\neq\chi\in H^1 (G_S, \Z/ \ell \Z )$, we call
$\sigma_u (\chi )=
< \chi_u , 1+{\ell}_u >$ the $u$-{\em signature} of $\chi$.\\

Let $g\in\Z$ so that $g\bmod p$ generates the multiplicative group
of $\F_p$.  Then the class of $g$ generates $\coKv^* /\coKv^{*\ell}
\cong \Z_{p}^* /\Z_{p}^{*\ell} \cong \Z/{\ell}\Z$. For $\chi\in H^1
(G_S, \Z/ \ell \Z )$, we call $\sigma_v (\chi ) = < \chi_v , g_v >
\neq 0$ the
$v$-{\em signature} of $\chi$.\\

We call the pair $(\sigma_u (\chi), \sigma_v (\chi))$ the {\em
signature} of $\chi$.  If we take $\chi'$ satisfying the conditions
we have put, then $\chi'=a\chi$ for some $a\in(\Z/{\ell}\Z)^*$, and
hence we don't change $\sigma_u (\chi) \sigma_v
(\chi)^{-1}\in\Z/{\ell}\Z$.  This last quantity only depends on
$K_S$ and we call it
 the {\em ramification signature} of $K_S$; it is nonzero.

 \subsection{DL and Signature Computation}
 \label{equivalence}
In this section we show that the discrete logarithm problem in the
multiplicative case is random polynomial time equivalent to
computing the signature of cyclic extensions with prescribed
ramification as described in Proposition~\ref{oneram}.\\
\ \\{\bf DL Problem}:  Suppose we are given $p$, $\ell$, $g$ and
$a$, where $p$ and $\ell$ are prime with $p\equiv 1 \pmod{{\ell}}$,
$g$ is a generator for the group $\F_p^*\{\ell\}$  of elements
killed by $\ell$, and $a\in\F_p^*\{\ell\}$. Then compute $ m \bmod
{\ell}$ such that $a = g^m$ in $\F_p$.\\

\ \\{\bf Signature Computation Problem}: Suppose we are given $K$,
$p$, $\ell$, $u$, $u'$, $v$, $\alpha$ and $g$, where
 $K=\Q(\sqrt{D})$ is a real quadratic field, $\ell$, $p$ are primes that split in
 $K$, $u$ and $u'$ are the two places of $K$ over $\ell$, $v$ is a place of $K$ over $p$, $\alpha$ is a
 unit of $K$, and $g$ is
a generator for $\F_p^*$ such that: (1)
 the class number of $K$ is not divisible by $\ell$,
(2) $\alpha^{l-1} \not\equiv 1 \pmod{P_u^2}$, $\alpha^{l-1}
\not\equiv 1 \pmod{P_{u'}^2}$, and (3) $\alpha^{\frac{p-1}{{\ell}}}
\not\equiv 1 \pmod{P_v}$. Then compute the ramification signature,
with respect to $1+\ell$ and $g$, of the cyclic extension of degree
$\ell$ over $K$ which is ramified at $u,v$ and unramified elsewhere.

\begin{theorem}
The problems DL and Signature Computation are random polynomial time
equivalent.
\end{theorem}

For the proof of the theorem, we first give a random polynomial time
reduction from DL to Signature Computation.  This part of the proof
depends on some heuristic assumption which will be made
clear below.\\

Let $a = g^m$ in $\F_p$ where $m$ is to be computed. If
$a^{\frac{p-1}{\ell}}=1$, then $m\equiv 0\pmod{\ell}$.  So suppose
$a^{\frac{p-1}{\ell}}\neq 1$.  We lift $a$ to some unit $\alpha$ of
a real quadratic field $K$ such that $\alpha\equiv a \pmod{v}$ for
some place $v$ of $K$ over $p$. This can be done as follows.
\begin{enumerate}
\item Compute $b\in\F_p$ such that $ab = 1$ in $\F_p$. \item
 $c \leftarrow 2^{-1} (a+b)$; $d\leftarrow 2^{-1} (a-b)$.
Note that $c^2 - d^2 =1$, and $a=c+d$.  We may assume $d\neq 0$
otherwise $a^2 = 1$ and $m=(p-1)/2$ or $p-1$. \item
 Treat $d$ as an integer.  Let $\gamma\in \bar{\Q}$ be such that $\gamma^2 = 1+d^2$.
\item Check if $1+d^2$ is a quadratic residue modulo ${\ell}$.
Otherwise substitute $d+rp$ for $d$ for random $r$ until the
condition is met.  This is to make sure that ${\ell}$ splits in $K$.
\item $\gamma^2 = 1+d^2 \equiv c^2 \pmod {p}$ implies $\gamma \equiv
c \pmod {v}$, and $\gamma \equiv -c \pmod {v'}$ where $v$ and $v'$
are the two places of $K$ over $p$. \item Let $\alpha = \gamma +d$.
Then $\alpha \equiv c+d \equiv a \pmod{v}$.  Note that the norm of
$\alpha$ is $d^2 - \gamma^2 = -1$, so $\alpha$ is a unit of $K$.
\end{enumerate}

We make the heuristic assumption that it is likely for $K$ to
satisfy the conditions in Proposition~\ref{oneram}. (Note that
condition (3) is satisfied since $\alpha \equiv a\pmod{v}$ and
$a^{\frac{p-1}{{\ell}}} \neq 1$.) We argue below that computing the
discrete logarithm $m$ where $a = g^m$ is reduced to solving the
Signature Computation problem on input $K$, $p$, $\ell$, $u$, $v$,
$\alpha$ and $g$, where $K=\Q(\gamma)$ with $\gamma^2 = 1+d^2$,
$\alpha = \gamma+d$, $u$ and $v$ are as constructed above. A simple
analysis shows that the expected time complexity in constructing
these objects is $O(\log^3 p )$.\\

For $\chi\in H^1 (G_S,\Z/{\ell}\Z)$ that is ramified at $u$ and $v$,
and unramified elsewhere, we have
\[ 0= < \chi_u , \alpha_u > + < \chi_v , \alpha_v > .  \]
Moreover since $\alpha^{\frac{p-1}{{\ell}}}\not\equiv 1\pmod{v}$,
$\alpha$ generates $\coKv^*/\coKv^{*\ell}$, so $< \chi_v , \alpha_v
> \neq 0$, and it
follows that $< \chi_u , \alpha_u > \neq 0$.\\

In general for a field $k$ and $a,b\in k^*$, we write $a\sim^l b$
if $a/b \in k^{*{\ell}}$.\\

We have $\alpha \sim^l g^m$ in $K_v$ since $\alpha\equiv a \equiv
g^m \pmod{v}$.  Hence $$< \chi_v , \alpha_v
> = < \chi_v , g^m_v > = m < \chi_v  , g_v  >.$$\\

Write $\alpha = \xi (1+y{\ell}) \pmod{{\ell}^2}$ with
$\xi^{{\ell}-1} = 1$ after identifying $\alpha$ with its isomorphic
image in $\Q_{\ell}$.   Then $\alpha \sim^{\ell} (1+{\ell})^y $, and

\[ 0= < \chi_u , \alpha_u >  = < \chi_u , (1+{\ell})^y_u >
= y < \chi_u , 1+{\ell}_u > .  \]

Hence we have
\[ < \chi_u , \alpha_u > + < \chi_v , \alpha_v >
= y < \chi_u , 1+{\ell}_u > + m < \chi_v , g_v > .  \]

So $y\sigma_u (\chi ) + m \sigma_v (\chi ) =0$.  From this we see
that if the ramification signature $\sigma_u (\chi)(\sigma_v
(\chi))^{-1}$ is determined then $m$ is determined.  The expected
time in this reduction is $O(\log^3 p)$. \\

Next we give a random polynomial time reduction from Signature
Computation on input $K$, $p$, $\ell$, $u$, $v$,
$\alpha$ and $g$, to DL on input $p$, $\ell$, $g$ and $a$ where $\alpha \equiv a \pmod{v}$.\\

Call the oracle to DL on input $p$, $\ell$, $g$ and $a$ to compute
$m$ such that $g^m = a \pmod{p}$.  Then $\alpha \equiv
g^m \pmod{v}$.\\

Write $\alpha = \xi (1+y{\ell}) \pmod{{\ell}^2}$ with
$\xi^{{\ell}-1} = 1$ after identifying $\alpha$ with its isomorphic
image in $\Q_{\ell}$.  Then $\alpha \sim^{\ell} (1+{\ell})^y $.
Again, $\xi \bmod{{\ell}^2}$ and hence $y$ can be computed
efficiently in time $O( | | \alpha | | \log\ell +\log^3 \ell)=O( | |
\alpha | | \log p
+\log^3 p)$.\\

For $\chi\in H^1 (G_S,\Z/{\ell}\Z)$ that is ramified at $u$ and $v$,
and unramified elsewhere, we have as before $< \chi_v , \alpha_v
> = < \chi_v  , g^m_v  >= m < \chi_v , g_v > ,$ and $< \chi_u , \alpha_u >  = < \chi_u , (1+{\ell})^y_u >
= y < \chi_u , 1+{\ell}_u > . $ Hence

\[ 0= < \chi_u , \alpha_u > + < \chi_v , \alpha_v >
= y < \chi_u , 1+{\ell}_u > + m < \chi_v  , g_v  >  \] from this we
can determine the signature $\sigma_u (\chi)(\sigma_v (\chi))^{-1}$.
The expected running time in this reduction is $O( | | \alpha | |
\log p +\log^3 p)$

\section{Signature Calculus for ECDL}
\subsection{Preliminaries}

In this section we will demonstrate the existence of principal
homogeneous spaces of order $\ell$ under elliptic curves over number
fields with prescribed ramification.  We begin by describing
$H^1(K_v,E)[\ell]$ in general terms when $E$ has good reduction at
$v$.

\begin{lemma}
\label{h1ellip} Let $K_v$ be a local field with finite residue field
$k$.  Let $E$ be an elliptic curve defined over $K_v$ with good
reduction.
\begin{enumerate}
\item Suppose the characteristic of $k$ is $\ell$.  Then $H^1(K_v ,E)[\ell]\cong\Z/\ell\Z$ if
$K_v\cong\Q_{\ell}$ and $\ell \nmid \#\tilde{E}(k)$.
\item Suppose
the characteristic of $k$ is not $\ell$. Then
\begin{enumerate}
\item $H^1 (K_v , E)[\ell]=0$ if $\ell \nmid \#\tilde{E}(k)$;
\item $H^1(K_v ,E)[\ell]\cong\Z/\ell\Z$ if $\ell \mid
\#\tilde{E}(k)$ but $\ell^2 \nmid \#\tilde{E}(k)$.
\end{enumerate}
\end{enumerate}
\end{lemma}

\ \\{\bf Proof}
%Note that for $\chi\in H^1 (K_v, E)$.  Then $\chi$
%is unramified iff $\chi = 0$.

Let $E_1(K_v)$ be the kernel of the reduction map from $E(K_v)$ to
$\tilde{E}(k)$.  From the commutative diagram
\[
\begin{array}{clclclcll}
0 & \rightarrow & E_1 (K_v)    & \rightarrow & E(K_v)          & \rightarrow & \tilde{E} (k) & \rightarrow & 0\\
  &             & \downarrow \ell &             & \downarrow \ell&             & \downarrow \ell &            &  \\
0 & \rightarrow & E_1 (K_v)    & \rightarrow & E(K_v)          &
\rightarrow & \tilde{E} (k) & \rightarrow & 0
\end{array}
\]
 and the snake lemma, we get the exact sequence
\[
\begin{array}{c}
0 \rightarrow E_1 (K_v)[\ell] \rightarrow E(K_v) [\ell] \rightarrow
\tilde{E}(k)[\ell] \rightarrow E_1 (K_v)/ \ell E_1(K_v)
\\
\rightarrow E(K_v)/{\ell}E(K_v) \rightarrow \tilde{E}(k)/ \ell
\tilde{E}(k) \rightarrow 0.
\end{array}
\]

If $\ell$ does not divide the order of $\tilde{E}(k)$, then
$\tilde{E} (k) [\ell ]$ and $\tilde{E}(k)/ \ell \tilde{E}(k)$ are
both 0.
Hence $ E(K_v)/{\ell}E(K_v) \cong E_1 (K_v)/{\ell}E_1 (K_v)$.\\

To prove (1) suppose the characteristic of $k$ is $\ell$. If $K_v
\cong \Q_{\ell}$ , then $E_1 (K_v )/\ell E_1 (K_v ) \cong \Z /\ell
\Z$. If moreover $| \tilde{E}(k) |$ is not divisible by $\ell$, then
$ E(K_v)/{\ell}E(K_v) \cong E_1 (K_v)/{\ell}E_1 (K_v) \cong \Z /\ell
\Z$, hence $H^1 ( K_v ,
E)[\ell] \cong \Z/\ell\Z$ by local duality.\\

To prove (2), suppose the characteristic of $k$ is not $\ell$. Then
$E_1 (K_v)/ \ell E_1(K_v)=0$, and it follows from the long exact
sequence that $ E(K_v)/{\ell}E(K_v) \cong \tilde{E}
(k)/{\ell}\tilde{E} (k)$. If $\ell$ does not divide the order of
$\tilde{E}(k)$, then $ E(K_v)/{\ell}E(K_v) \cong \tilde{E}
(k)/{\ell}\tilde{E} (k)=0$, and by local duality, $H^1 (K_v, E)
[\ell ] =0$. This proves 2(a). If $| \tilde{E} (k ) |$ is divisible
by $\ell$ but not $\ell^2$, then $\tilde{E}(k)/ \ell
\tilde{E}(k)\cong \Z/\ell \Z$.  Hence $ E(K_v)/{\ell}E(K_v) \cong
\tilde{E} (k)/{\ell}\tilde{E} (k) \cong \Z/\ell \Z$, and by local
duality, $H^1 (K_v, E) [\ell ] = \Z /\ell \Z$. Thus 2(b) is proved.

\vspace{.2in}
\subsubsection{Ranks of Quadratic Twists of Elliptic Curves over
$\Q$}

Let $E$ be an elliptic curve over $\Q$ and fix a Weierstrass
equation for $E$:

$$y^2=x^3+ax+b.$$
Let $K=\Q(\sqrt{d})$ be a quadratic extension of $\Q$ and let $E_d$
be the quadratic twist of $E$ given by the equation

$$dy^2=x^3+ax+b.$$

Let $G$ be the Galois group of $K$ over $\Q$ and $\sigma$ a
generator of $G$. Denote by $V$ the group $E(K)\otimes_\Z\Q$, by
$V^+$ the fixed space by $\sigma$, and by $V^{-}$ the subspace of
$V$ where $\sigma$ acts by $-1$.  Now

$$V=V^{+}\oplus V^-,$$

$V^+=E(\Q)\otimes_{\Z}\Q$, and we see easily that
$V^-=E_d(\Q)\otimes_{\Z}\Q$, via the isomorphism sending a point
$(x,y)$ in $E_d(\Q)$ to $(x,\sqrt{d}y)$ in $V^-$.\\

In the algorithm in \S~\ref{sgecdl} below, it will help to have a
lifting $E/\Q$ of our original elliptic curve $\tilde{E}/\F_p$ such
that both $E(\Q)$ and $E_d(\Q)$ are of rank one. Standard
conjectures about the behavior of the rank of the Mordell-Weil group
of an elliptic curve predict that it should be quite possible to
find such a situation. For example, a conjecture of Goldfeld [G]
says that the rank of a quadratic twist $E_d$ of an elliptic curve
$E$ over $\Q$ should be on average as small as the sign of the
functional equation of its $L$-function would allow, i.e. either 0
or 1, depending on whether this sign is +1 or -1.  In fact, assuming
the Riemann hypothesis for all of the curves $E_d$, Heath-Brown [HB]
has proved that at least $1/4$ of all the $E_d$ with the sign in the
functional equation of the L-function being +1 will have rank 0 and
at least 3/4 of all the $E_d$ with the sign being -1 will have rank
1 (see [H-B], Theorem 4).  In the algorithm, we will first lift
$\tilde{E}/\F_p$ to $E/\Q$ that has rank at least one by
construction, and we will make the heuristic assumption that $E(\Q)$
is of rank exactly one and therefore the sign of the functional
equation is -1 (see [BMSW], \S 3 for why this is considered to be
reasonable). Using ([RS], Theorem 7.2), Heath-Brown's result just
mentioned, and taking sufficiently many random $d$, we can
heuristically arrange for the sign of the functional equation of
$E_d$ to be equal to -1 and for $E_d(\Q)$ to have rank 1. When we
make the heuristic assumption in \S~\ref{sgecdl} below that the rank
of $E(K)$ is exactly two, we shall mean this.

\vspace{.2in}

\subsubsection{The Group $E(K_v)/\ell$ at Bad Reduction Primes $v$ of
$E$}

Let $\tilde{E}$ be an elliptic curve over $\F_p$, $p\geq 5$, given
in Weierstrass form by an affine equation

$$y^2=x^3+\tilde{a}x+\tilde{b}.$$ In the algorithm below, we will want to lift $\tilde{E}$ to an
elliptic curve $E$ over $\Q$ with Weierstrass equation

$$y^2=x^3+ax+b,$$
having good reduction at $\ell$ and such that at primes $v$ of bad
reduction, $E(K_v)/\ell=0$.  We give a heuristic here about why this
should be possible.  In our lifting in the algorithm presented in
\S~\ref{sgecdl}, $| a |$ is at most $p^2$ and $| b |$ is at most
$p^4$, so the discriminant $\Delta$ of a minimal Weierstrass
equation for $E$ is of order at most $p^8$.  At a prime $v$ of split
multiplicative reduction prime, the group of connected components of
the N\'eron model of $E$ over the ring of integers of $\Q_v$ will be
of order the power of $v$ in the discriminant.  Since $\ell$ is of
the same order as $p$, this power is very unlikely to be divisible
by $\ell$. At other primes of bad reduction, the group of connected
components is of order at most 4 (see [Si1], Ch. VII, Theorem 6.1).
Thus the order of the group of components is unlikely to be
divisible by $\ell$.  We claim that this implies that for any bad
reduction place $v$ of $E$, $E(K_v)/\ell=0$.  To see this, recall
that $E(K_v)$ has a filtration:

$$E(K_v)\supseteq E_0(K_v)\supseteq E_1(K_v),$$
where $E_0(K_v)$ is the group of points specializing to points of
the smooth locus $\tilde{E'}^0$ of the special fibre $\tilde{E'}$ of
the minimal regular proper model $\mathcal{E}$ of $E$ over the ring
of integers $R$ of $K_v$ and $E_1(K_v)$ is the kernel of the
reduction map

$$E_0(K_v)\to \tilde{E'}^0(\F).$$  Now $E(K_v)/E_0(K_v)$ is the group of connected components of
the special fibre of the N\'eron model of $E$ over $R$, and $E_1(K)$
is a pro-$v$-group, where $v$ is the residue characteristic of
$K_v$. $E_0(K_v)/E_1(K_v)$ is the group of points on the connected
component of identity of the special fibre of the N\'eron model.
This last group is isomorphic to either the additive group or the
multiplicative group of the residue field, $\F_v$.  Because $E$ has
good reduction at $\ell$, $v\neq \ell$, $\ell$ is large and $v$ is
relatively small compared to $\ell$, it is unlikely that $\ell$ will
divide $v-1$. Thus it is likely that $E(K_v)/\ell=0$ and we shall
use this heuristic in the algorithm in \S~\ref{sgecdl} below.

\subsection{Principal Homogeneous Spaces Ramified over $p$ and
$\ell$} Throughout this section, let $p,{\ell}$ be odd, rational
primes. Let $K/\Q$ be a real quadratic extension,
$X=\mbox{Spec}(\co_K)$ and $\Sigma$ be the set of all places at
which $E$ has bad reduction, together with all the archimedean
places. Let ${\mathcal E}$ be a smooth proper model of $E$ over the
open subset $U=X-\Sigma$.  If $S$ is any set of places of $K$
containing $\Sigma$, denote by $U_S$ the open set $X-S$.  We denote
by $\hbox{\cyr Sh}(E)$ the Shafarevich-Tate group of everywhere
locally trivial principal homogeneous spaces under $E$ over $K$.

\begin{proposition}
\label{seq} Let $S$ be a finite set of places of $K$ containing all
bad reduction places of $E$ and the places above $\ell$.  Then if
$\hbox{\cyr Sh}(E)\{\ell\}=0$, we have the exact sequence:

$$E(K)/\ell\to \prod_{v\in S}E(K_v)/\ell\to (H^1(U_S,{\mathcal
E})[\ell])^*\to 0.$$
\end{proposition}

Proof:  Consider the Cassels-Tate exact sequence

$$(**)\: E(K)^{(\ell)}\to \prod_{v\in S}E(K_v)^{(\ell)}\to H^1(U_S,\cE)\{\ell\}^*\to \hbox{\cyr
Sh}(E)\{\ell\}\to 0.$$

\begin{lemma}Let $B$ be a torsion abelian group such that $B[\ell^n]$ and $B/\ell^nB$ are
finite groups.  Then we have

$$B[\ell]^*\cong B^*/\ell B^*$$

and

$$B\{\ell\}^*\cong B^{*(\ell)}$$
\end{lemma}

Proof:

Let $n$ be a positive integer and consider the tautological exact
sequence:

$$0\to B[\ell^n]\to B\stackrel{\ell^n}{\to} B\to B/\ell^n B\to 0.$$

Since the functor * (see \S 1 for notation) is exact on the category
of locally compact abelian groups, we get the exact sequence:

$$0\to (B/\ell^n B)^*\to B^*\stackrel{\ell^n}{\to}B^*\to B[\ell^n]^*\to
0.$$ We then get the first conclusion of the lemma by taking $n=1$
and the second by passing to the inverse limit over $n$.  This
completes the proof of the lemma.\\

The proposition then follows from the lemma, the assumption that
$\hbox{\cyr Sh}(E)\{\ell\}=0$, and the
Cassels-Tate sequence above by reducing the terms mod $\ell$.\\

For the remainder of this section we assume that $p$ and ${\ell}$
split in $K$, and that $E$ has good reduction at $p$ and $\ell$,
with $\#\tilde{E} (\F_p)=\ell$ and $\ell\neq\#\tilde{E}(\F_{\ell})$.
Moreover, because we assume that $\ell$ is sufficiently large, a
theorem of Kamienny [Ka] ensures that $E(L)[\ell]$ is trivial for
all quadratic extensions $L$ over $\Q$. Finally, we assume that
$E(K_v)/\ell=0$ for all bad reduction places $v$ of $E$ (see \S 1.3
for why this is reasonable, heuristically).

\begin{proposition}
\label{eram} Let $S$ be a finite set of places of $K$ containing all
bad reduction places of $E$ and the places above $\ell$.  $S$ may or
may not contain places above $p$, but assume that it contains no
good reduction places that do not divide $\ell$ or $p$. Suppose
\begin{enumerate}
\item $\hbox{\cyr Sh}(E)\{\ell\}=0$;
\item the map $E(K)/\ell\rightarrow E(K_u)/\ell\oplus
E(K_{u'})/\ell$ is an isomorphism, where $u$ and $u'$ are the two
places of $K$ over $\ell$.
\end{enumerate}
Then the $\F_{\ell}$-dimension of  $H^1 (U, {\mathcal E})[\ell]$
equals $n(S)-2$ where $n(S)$ is the number of finite places in
$S-\Sigma$.
\end{proposition}

Proof: Since $\hbox{\cyr Sh}(E)\{\ell\}=0$, we have the exact
sequence  $$E(K)/\ell\to \prod_{v\in S}E(K_v)/\ell\to
(H^1(U_S,{\mathcal E})[\ell])^*\to 0$$ by Proposition~\ref{seq}. The
middle group in the sequence $\prod_{v\in S}E(K_v)/\ell$ is
isomorphic to the direct sum of $n(S)$ copies of $\Z/{\ell}\Z$ by
Lemma~\ref{h1ellip}. Since the map
$$E(K)/\ell\rightarrow
E(K_u)/\ell\oplus E(K_{u'})/\ell\cong \Z/\ell\Z \oplus \Z/\ell\Z$$
is an isomorphism, it follows that the image of the map
$$E(K)/\ell\to \prod_{v\in S}E(K_v)/\ell$$ is isomorphic to $\Z/{\ell}\Z\oplus
\Z/\ell\Z$. Hence the $\F_{\ell}$-dimension of $H^1 (U_S, {\mathcal
E})[\ell]$ equals $n(S)-2$.

\begin{proposition}
\label{eoneram} Let $S$ be the set consisting of all bad reduction
places of $E$, together with the two places $u$ and $u'$ over
$\ell$, and one place $v$ over $p$. Suppose
\begin{enumerate}
\item $\hbox{\cyr Sh}(E)\{\ell\}=0$;
\item the map $E(K)/\ell\rightarrow E(K_u)/\ell\oplus
E(K_{u'})/\ell$ is an isomorphism.
\end{enumerate}
Then the $\F_{\ell}$-dimension of  $H^1 (U, {\mathcal E})[\ell]$ is
one.  Moreover, every nontrivial element of $H^1 (U, {\mathcal
E})[\ell]$ is ramified at $v$.
\end{proposition}
\ \\{\bf Proof} Suppose $u,u'$ are the places over $\ell$, $v,v'$
the places over $p$.  Let $R=\Sigma\cup\{u,u'\}$ and
$T=\Sigma\cup\{u,u',v\}$.  Then from Proposition~\ref{eram} we know
that $H^1 (U_R ,\E )[\ell]$ has dimension zero and $H^1 (U_T ,\E
)[\ell]$ has dimension one. So there exists $\chi\in H^1 (U_T ,\E
)[\ell] - H^1 (U_R ,\E )[\ell]$ and  $\chi$ must be ramified at
$v$.  This completes the proof of the proposition.\\

We remark that in the proposition above the assumption that the map
$E(K)/\ell\rightarrow E(K_u)/\ell\oplus E(K_{u'})/\ell$ is an
isomorphism can be replaced by the assumption that the image of
$E(K)/\ell$ in $E(K_u)/\ell\oplus E(K_{u'})/\ell$ and in
$E(K_v)/\ell \oplus E(K_u)/\ell\oplus E(K_{u'})/\ell$ are both of
$\F_{\ell}$-dimension
two.\\

For $w=u,u',v$, let $R_w\in E(K_w)-\ell E(K_w)$, so that the class
of $R_w$ generates $E(K_w)/\ell$. For $\chi\in H^1 (U_S, {\mathcal
E})[\ell]$ and $w$ a place of $K$, we call $a_w = < \chi_w, R_w >$
the $w$-$signature$ of $\chi$.  We call $(a_u, a_{u'}, a_v)$ the
signature of $\chi$ with respect to $R_u$, $R_{u'}$ and $R_v$.
Proposition~\ref{eoneram} implies that $< \chi_v , R_v
> \neq 0$ for any nontrivial $\chi\in H^1 (U_S, {\mathcal
E})[\ell]$, and $( \frac{< \chi_u , R_u
>}{ < \chi_v , R_v > },\frac{< \chi_{u'} , R_{u'}
>}{ < \chi_v , R_v > })$ is the same for all such $\chi$.  We call
this pair the {\em signature} of $H^1 (U, {\mathcal E})[\ell]$
with respect to $R_u$, $R_{u'}$ and $R_v$.\\

Since the pairing between $H^1 (K_v, E)[\ell]$ and $E(K_v)/\ell
E(K_v)$ is perfect, both being isomorphic to $\Z_{\ell}/\Z_{\ell}$,
there is a unique $\psi_v\in H^1 (K_v,E)[\ell]$ such that $< \psi_v,
R_v
> =1$. Similarly, there is a unique $\psi_{u}\in H^1
(K_{u},E)[\ell]$ such that $< \psi_{u}, R_{u}
> =1$, and a unique
$\psi_{u'}\in H^1 (K_{u'},E)[\ell]$ such that $< \psi_{u'}, R_{u'}
> =1$.
Let $\chi\in H^1 (U_S, {\mathcal E})[\ell]$. Suppose $\chi_v = a_v
\psi_v$, $\chi_u = a_u\psi_u$ and $\chi_{u'}=a_{u'}\psi_{u'}$. Then
$< \chi_v, R_v > = a_v$, $< \chi_u , R_u > = a_u$ and $< \chi_{u'} ,
R_{u'} > = a_{u'}$. So $a_u$, $a_{u'}$ and $a_v$ constitute the
signature for $\chi$ with respect to $R_u$, $R_{u'}$ and $R_v$. Thus
the signature $(a_u, a_{u'}, a_v)$ succinctly represents the
localization of $\chi$ at the ramified places. These localizations
in turn determine $\chi$ uniquely, since the Shafarevich-Tate group
is assumed to have trivial $\ell$-part. Therefore, the signature of
$\chi$ can be regarded as a succinct representation of $\chi$  (by
determining its localization at $u$, $u'$ and $v$ as $\chi_u =
a_u\psi_u$, $\chi_{u'}=a_{u'}\psi_{u'}$, and $\chi_v = a_v \psi_v$).
We note that this representation requires only $O(\log\ell)$ bits
whereas an explicit description of
$\chi$ may require $\Omega(\ell)$ bits.\\

Suppose in addition to the map $E(K)/\ell\rightarrow
E(K_u)/\ell\oplus E(K_{u'})/\ell$ being an isomorphism, we assume
that the map $E(K)/\ell\rightarrow E(K_v)/\ell$ is nontrivial.  In
this case we may obtain $R_w$'s as follows. Let $Q,R\in E(K)$ so
that their classes generate $E(K)/\ell$. Suppose without loss of
generality that the class of $Q$ is nontrivial in $E(K_u)/\ell$ and
the class of $R$ is nontrivial in $E(K_{u'})/\ell$. As
$E(K)/\ell\rightarrow E(K_v)/\ell$ is nontrivial, the class of
either $Q$ or $R$ is nontrivial in $E(K_v)/\ell$. Suppose without
loss of generality the class of $Q$ is nontrivial in $E(K_v)/\ell$.
Then we may take $R_v = Q$, $R_u = Q$ and $R_{u'}=R$.\\

\subsection{ECDL and Signature Computation}
\label{sgecdl}
 In this section we show that the elliptic curve discrete logarithm problem
is random polynomial time equivalent to computing the signature of
homogeneous spaces with prescribed
ramification as described in Proposition~\ref{eoneram}.\\

\ \\{\bf ECDL}:  Given $p$, $\ell$, $\tilde{E}$, $\tilde{Q}$ and
$\tilde{R}$, where $p$ and  $\ell$ are prime,  $\tilde{E}$ is an
elliptic curve defined over $\F_p$ with $\#\tilde{E}(\F_p)=\ell$,
and non-zero points $\tilde{Q},\tilde{R}\in\tilde{E}(\F_p)$, to
determine $m$
so that $\tilde{R}=m \tilde{Q}$.\\

\ \\ {\bf Homogeneous Space Signature Computation}: Suppose we are
given $p$, $\ell$, $K$, $E$, $v$, $Q$, $R$,  where $p$ and $\ell$
are prime, $K$ is a quadratic field where $p$ and  $\ell$ both
split, $E$ is an elliptic curve defined over $K$ with $\hbox{\cyr
Sh}(E)\{\ell\}=0$ and the discriminant of $E$ being prime to $\ell$,
$v$ is a place of $K$ over $p$, $Q$ and $R\in E(K)$ such that
$Q\not\equiv 0 \pmod{\ell E(K_v)}$ and the images of $R$ and $Q$ in
$E(K_u)/\ell\oplus E(K_{u'})/\ell$ are independent, where $u$ and
$u'$ are the two places of $K$ over $\ell$.  Then compute the
signature of $H^1 (U_S, {\mathcal E})[\ell]$ with respect to $\rho_v
= Q$, $\rho_u = Q$ and $\rho_{u'}=R$, where $S$ is the set
consisting of $u,u',v$ and all places of bad reduction of $E$. (Note
that $\rho_w$ generates $E(K_w)/\ell E(K_w)$ for $w = u, u' , v$.)

\begin{theorem}
The problems ECDL and Homogeneous Space Signature Computation are
random polynomial time equivalent.
\end{theorem}

For the proof of the theorem, we first give a random polynomial time
reduction from ECDL to Homogeneous Space Signature Computation. This
part of the proof depends on some heuristic assumptions which will
be made
clear below.\\

Given $\tilde{E}/\F_p$ where $\tilde{E}(\F_p) [\ell]= < \tilde{Q}
>$, and
$\tilde{R}$, we are to compute $m$ so that $\tilde{R} = m\tilde{Q}$.
Steps 1-3 of the reduction construct an instance $p$, $\ell$, $K$,
$E$, $v$, $Q$, $R$ of the Homogeneous
Space Signature Computation problem.\\

1. Construct $E/\Q$ with $Q\in E(\Q)$ such that $\tilde{Q}= Q
\bmod{p}$. This can be done as follows. Suppose $\tilde{E}$ is
specified by an affine equation $y^2=x^3+\bar{a}x+\bar{b}$ where
$\bar{a}=a \bmod{p}$, $\bar{b}=b\bmod{p}$ with $0 \le a,b < p$ and
$\tilde{Q}=(u\bmod{p},v\bmod{p})$ with $0 < u,v < p$. Choose a
random integer $r$, $0\le r < p$,  and let $Q=(u,v+rp)$. Let $b_r =
(v+rp)^2 - (u^3+au)$.  Then $Q\in E_r (\Q)$ where $E_r$ is the
elliptic curve with the affine equation $y^2=x^3+ax+b_r$. Set
$E=E_r$. The point $Q$ cannot be torsion for otherwise it would have
to be in $E(\Q)[\ell]$, which has no non-zero point since $\ell$ is
big. The height of $Q$ is far smaller than that of a point in $\ell
E(\Q)$, so $Q$ is not in $\ell E(\Q)$. Since $\tilde{E}(\F_p)
[\ell]\cong \Z/\ell\Z$, $E(\Q_p)/\ell \cong \tilde{E} (\F_p )/\ell
\cong \Z/\ell\Z$ and the class of $Q$ generates $E(\Q_p)/\ell$.\\

2. Check that $E$ has good reduction at $\ell$ and that $|
\tilde{E}(\F_{\ell} ) |$ is not divisible by $\ell$. Otherwise, go
back to 1. to find a different $E$.\\

3.  Lift $\tilde{R}$ to $R\in E(K)$ where $K/\Q$ is a quadratic
extension in which $p$ and $\ell$ both split. This can be done as
follows.  Suppose $E$ is defined by the affine equation
$y^2=x^3+ax+c$.  Suppose $\tilde{R}=(\mu\bmod{p}, \nu\bmod{p})$ with
$0 < \mu, \nu < p$. Choose a random positive integer $r < p$. Set
$\mu_r = \mu+rp$. Let $\beta$ be a root of $y^2 = \mu_r^3+ a\mu_r +
c$.  Then $(\mu_r,\beta)$ is a lift of $\tilde{R}$ in $E(K)$ where
$K=\Q(\beta)$.  By construction $p$ splits in $K$,
$$E(K_v)/\ell \cong E(\Q_p)/\ell \cong \tilde{E} (\F_p )/\ell \cong
\Z/\ell\Z$$ and $R - mQ \in \ell E(K_v)$.  Check that $\ell$ splits
in $K$ and that the images of $R$ and $Q$ in $E(K_u)/\ell\oplus
E(K_{u'})/\ell$ are independent;  otherwise repeat the above steps
with a different $r$ until a suitable $K$ is found. Say the class of
$Q$ is nontrivial in $E(K_u)/\ell$ and the class of $R$ is
nontrivial
in $E(K_{u'})/\ell$.\\

4.  Call the oracle for the Homogeneous Space Signature Computation
on input $p$, $\ell$, $E$, $K$, $Q$, $R$, $v$ to compute the
signature $(\alpha,\beta)$ of $H^1 (U_S, {\mathcal E})[\ell]$ with
respect to $\rho_v = Q$, $\rho_u = Q$ and $\rho_{u'}=R$ (where $S$
is the set consisting of $u,u',v$ and all places of bad reduction of
$E$). Then for all nontrivial $\chi\in H^1 (U_S, {\mathcal
E})[\ell]$, $\alpha = \frac{ < \chi_u , Q_u >}{ < \chi_v , Q_v
>}$ and $\beta=\frac{ < \chi_{u'} , R_{u'} >}{ < \chi_v , Q_v > }$.\\

5. Identify $K_{u}$ with $\Q_{\ell}$ and compute $n$ so that $R
\equiv n Q \pmod{\ell E(K_{u} )}$ as follows. Compute $d = |
\tilde{E} (\F_{\ell} ) |$. Observe that $dQ$ and $dR$ are both in
$E_1 (\Q_{\ell})$. Compute $n$ such that $n (dQ) \equiv (dR)
\pmod{\ell}$ in $E_1 (\Q_{\ell} )$. Then $d (nQ-R) = \ell Z$ for
some $Z\in E_1 (\Q_{\ell} )$. Since $d$ is not divisible by $\ell$,
$d^{-1}\in \Z_{\ell}$, so $nQ - R = d^{-1} \ell
Z = \ell (d^{-1} Z) \in \ell E(\Q_{\ell})$.\\

6.  Now
\begin{eqnarray*}
0 & = & \sum_{w\in\{v,u,u'\}} < \chi_w , R_w > \\
  & = & m < \chi_v  , Q_v  >+  n < \chi_u , Q_u > + < \chi_{u'}  , R_{u'}  >.
\end{eqnarray*}
From this we get $m+n\alpha+\beta \equiv 0\pmod{\ell}$.  Hence $m$
can be determined.\\

We make the heuristic assumption that it is likely for $E$ and $K$
to satisfy the conditions in Proposition~\ref{eram}.  Note that by
construction $E(\Q)$ is of rank at least one.  The points $Q$ and
$R$ are likely to be integrally independent in $E(K)$ as they both
have small height by construction.  So $E(K)$ is likely to be of
rank at least two and we make the heuristic assumption that with
nontrivial probability its rank is exactly two. Moreover, since
$Q\in E(\Q)$ and $R\in E(K)-E(\Q)$, the images of $Q$ and $R$ are
likely to be independent in $E(K_u)/\ell \oplus E(K_{u'})/\ell$,
heuristically speaking.  The expected running time of this reduction
is dominated by Step 2 where the number of rational points on the
reduction of $E$ mod $\ell$ is
counted.  The running time of that step is $O(\log^8 \ell)$ [Sc], hence it is $O(\log^8 p )$.\\

Next we give a random polynomial time reduction from Homogeneous
Space Signature Computation with input $p$, $\ell$, $E$, $K$, $Q$,
$R$, $v$ to ECDL with input $p$, $\ell$, $\tilde{E}$, $\tilde{Q}$,
$\tilde{R}$, where $\tilde{E}$ is the reduction of $E$ mod $v$,
$\tilde{Q}$ (resp.
$\tilde{R}$) is the reduction of $Q$ (resp. $R$) mod $v$.\\

For any nontrivial $\chi\in H^1 (K,E)[\ell]$ that is unramified away
from $u$, $u'$ and $v$, we have
\begin{eqnarray*}
< \chi_v  , Q_v  >+  < \chi_u , Q_u > + < \chi_{u'} , Q_{u'} > & = & 0,\\
< \chi_v , R_v > +  < \chi_u , R_u > + < \chi_{u'} , R_{u'} > & = &
0.
\end{eqnarray*}

Suppose $Q=a_w \rho_w \pmod{\ell E(K_w)}$ and $R=b_w \rho_w
\pmod{\ell E(K_w)}$ for $w=u, u', v$.  Note that from
Lemma~\ref{h1ellip}, $a_v$ and $b_v$ can be computed by solving ECDL
on the reduction of $E$ modulo $v$.  On the other hand $a_w$, $b_w$
for $w=u, u'$, can be computed in a manner as described in Step 5
above.

Then we get
\begin{eqnarray*}
a_v< \chi_v , \rho_v  > +  a_u < \chi_u , \rho_u > +  a_{u'} < \chi_{u'} , \rho_{u'} > & = & 0,\\
b_v< \chi_v  , \rho_v  >+  b_u < \chi_u  , \rho_u >+  b_{u'} <
\chi_{u'} , \rho_{u'} > & = & 0
\end{eqnarray*}

Condition (2) of Proposition~\ref{eram} implies that the two
relations above are linearly independent. From these we can compute
the the signature of $\chi$; that is $( \frac{< \chi_u , \rho_u
>}{ < \chi_v , \rho_v > },\frac{< \chi_{u'} , \rho_{u'}
>}{ < \chi_v  , \rho_v >})$.  The expected running time of this reduction
can be shown to be $O(\log^4 p)+O(M\log p )$ where $M$ is the
maximum of the lengths of $R$, $Q$ and $D$.

\section{Feasibility of Index Calculus}
We will derive an index calculus method for the signature
computation problem of Dirichlet characters.  We will discuss why a
similar method cannot work for principal homogeneous spaces.

\subsection{Index Calculus for Signature Computation of Dirichlet
Characters}Suppose we are given a real quadratic field $K$, primes
$\ell$, $p$, places $u,v$ satisfying the conditions in
Proposition~\ref{oneram}.  Let $K=\Q(\alpha)$ with $\alpha^2\in\Z_{>
0}$. To compute the signature of $\chi\in H^1 (K,\Z/\ell\Z)$ that is
ramified precisely at $u$ and $v$, we generate random algebraic
integers $\beta=r\alpha+s$ with $r,s\in\Z$ so that $r\alpha+s\equiv
g \pmod{v}$ and $\beta\sim (1+\ell)^a$ at $u$ for some $a$.  Now
suppose the norm of $\beta$ is $B$-smooth for some integer $B$. Then

\[0=\sum_w < \chi_w ,\beta_w  >= < \chi_v  , g_v  >+ a < \chi_u, 1+\ell_u >
+ \sum_w e_w < \chi_w , \pi_w > ,\]
where $w$ in the last sum ranges
over all places of $K$ of norm less than $B$, $\pi_w$ is a local
parameter at $w$, and $e_w$ is the valuation of $\beta$ under $w$.
Hence we have obtained a $\Z/\ell\Z$-linear relation on $(< \chi_v ,
g_v >)^{-1}< \chi_u, 1+\ell_u
> $, and  $(< \chi_v , g_v >)^{-1}< \chi_w, \pi_w
>$.  With $O(B)$ relations we can solve for all these unknowns, in
particular the signature $(< \chi_v , g_v >)^{-1}< \chi_u, 1+\ell_u
> $.
\subsection{The Elliptic Curve Case}
We see that one important reason why index calculus is viable in the
multiplicative case is due to the fact that locally unramified
Dirichlet characters can be paired nontrivially with non-units. For
the elliptic curve case, pairing a principal homogeneous space
$\chi$ and a global point $\alpha$ yields similarly a relation:
\[ 0=\sum_v < \chi_v , \alpha_v > . \]
However from Lemma~\ref{h1ellip} we see that in the sum above we
have nontrivial contribution from a place $v\nmid \ell$ (and where
$E$ has good reduction) only if $\ell$ divides $\#\tilde{E}(\F_v)$.
Since $\#\tilde{E}(\F_v)$ is of the order $\#\F_v$, which is the
norm of $v$, we see that the finite places of good reduction that
are involved in the sum are all of large norm. As for the bad
reduction places, the heuristic assumption that we discussed just
before Proposition~\ref{eram} implies that these will not play any
role in this sum, since it will be likely that $E(K_v)/\ell=0$ for
such places $v$, because $v$ is of small norm. This explains why the
index calculus method is lacking in the case of the elliptic curve
discrete logarithm problem.

\section{Characterization of ramification signature}
\label{signature}

Let $K, \ell, p, u, v, S$ be as in Proposition~\ref{oneram}.\\

Let $g\in\Z$ so that $g\bmod p$ generates the multiplicative group
of $\F_p$.  Let $w$ be the place of $K(\mu_{\ell})$ over $v$ such
that $g^{\frac{p-1}{{\ell}}} \equiv \zeta \pmod{w}$.\\

Let $M=K_S$ be the cyclic extension corresponding to $H^1 (G_S,
\Z/\ell\Z )$.  Suppose $\chi\in H^1 (G_S, \Z/\ell\Z )$ is
nontrivial.  Then $\chi$ corresponds to some $A\in K(\mu_{\ell})$
through $H^1 (K(\mu_{\ell}), \Z/{\ell}\Z)\cong H^1 (K(\mu_{\ell}),
\mu_{\ell}) \cong K(\mu_{\ell})^* /K(\mu_{\ell})^{*\ell}$, such that
$M(\mu_{\ell}) = K(\mu_{\ell}) (A^{\frac{1}{\ell}})$, and for all
$\sigma$ in the absolute Galois group of $K$, $\chi(\sigma)=i$ iff
$\sigma(A^{\frac{1}{{\ell}}}) / A^{\frac{1}{{\ell}}}=\zeta^i$.\\

The following proposition provides a concrete characterization of
the signature of $\chi$.

\begin{proposition}
\label{char} If we identify $K(\mu_{\ell})_w$ with $\Q_p$ and $K_u$
with $\Q_{\ell}$, then $A \sim^{\ell} p^{m}$ in $\Q_p^{ur}$ where
$m=\sigma_v (\chi)=<\chi_v, g_v >$, and $A\sim^{\ell} \zeta^n$ in
$\Q_{\ell}(\mu_{\ell})^{ur}$ where $n=\sigma_u (\chi)= < \chi_u ,
1+{\ell}_u >$.
\end{proposition}

The rest of this section is devoted to the proof of this
proposition.  We set some notation first. For any local field $L$,
let $L^{ur}$ denote the maximal unramified extension over $L$. For
any place $\nu$ of a number field $K$, let $\theta_{\nu}$ denote the
local Artin map, $$\theta_{\nu}: K_{\nu}^* \rightarrow G_{\nu}^{ab}
,$$ where $G_{\nu}^{ab}$ denotes the Galois group of the maximal
abelian extension of $K_{\nu}$. For $a,b\in K(\mu_{\ell})$ and $\nu$
a prime of $K(\mu_{\ell})$, we have $$\alpha^{\theta_{\nu} (b)} =
(a,b)_{\nu} \alpha$$ where $\alpha^l = a$, and $(a,b)_{\nu}$ denotes
the local norm residue symbol (see p.
351 of [CF]).\\

\begin{lemma}
$\sigma_u (\chi)= < \chi_u , 1+{\ell}_u > = \chi_u (\theta_u (1+\ell
))$ and $\sigma_v (\chi)=<\chi_v , g_v > = \chi_v(\theta_v (g))$
\end{lemma}
\ \\{\bf Proof} This follows directly from [S1] Chapter XIV
Proposition 3.

\ \\{\bf Proof of Proposition~\ref{char}}  Suppose $v'$ is a place
of $K(\mu_{\ell})$ such that $v' | v$. Then $d <\chi_v , b_v > =
<\chi_{v'} , b_{v'} >$ where $d= [ K(\mu_{\ell})_{v'} : K_v ]$ (see
[S], Proposition 7 of Ch. XIII). Moreover $<\chi_{v'} , b_{v'}
>=\chi_{v'}(\theta_{v'} (b))=i$ iff $(A,b)_{v'}=\zeta^i$.
Identifying $i$ with $\zeta^i$, we may write
\[ d <\chi_v , b_v > = <\chi_{v'} , b_{v'} > = (A,b)_{v'} \]\\

We analyze the situation at $p$ and $\ell$ separately.\\

(I) At $p$: $\Q_p^* / {\ell} = \mu_{{\ell}} \times < p > /{\ell} $.
So under the identification of $K(\mu_{\ell})_w$ with $\Q_p$, $A = u
p^{w(A)}$ where $u^{\ell} = 1$, and $e < {\ell}$.  Since $\Q_p
(u^{\frac{1}{{\ell}}})/\Q_p$ is unramified, $A \sim^{\ell}
p^{w(A)}$ in $\Q_p^{ur}$.\\

Let $\chi\in H^1 (K,\Z/{\ell}\Z )$
\[ < \chi_w, g_w > = ( A , g )_w = -(g,A)_w \]
\[ (g,A)_w = i \mbox{ iff } \zeta^i= \left( \frac{g}{w}\right)^{w(A)}\]
\begin{eqnarray*}
\left( \frac{g}{w}\right) & \equiv & g^{\frac{Nw - 1}{{\ell}}} \pmod{P_w}\\
                          & \equiv & g^{\frac{p-1}{{\ell}}} \pmod{P_w}\\
                           & \equiv & \zeta \pmod{P_w}
\end{eqnarray*}

Therefore, $(g,A)_w = w(A)$. Consequently,
$$< \chi_v , g_v > = < \chi_w , g_w > = -(g,A)_w = -w(A).$$\\

(II) At ${\ell}$:  Denote by $u'$ the place of $K(\mu_{\ell})$ over
$u$.  We have
\[({\ell}-1) < \chi_u , 1+{\ell}_u > = < \chi_{u'} , 1+{\ell}_{u'} > = (A,1+{\ell})_{u'}.\]
We verify below that $(A, 1+{\ell})_{u'} = n$.  Then we can conclude
that
$$\sigma_u (\chi ) = < \chi_u , 1+{\ell}_u > = -n.$$\\

There is a ramified extension of degree ${\ell}$ over $\Q_{\ell}$,
namely, the subextension $M_1$ of $\Q_{\ell} (\zeta_{{\ell}^2})$ of
degree ${\ell}$ over $\Q_{\ell}$.  Let $\psi$ be the ramified
character in $H^1 (\Q_{\ell} , \Z/{\ell}\Z )$ whose retriction to
$H^1 (\Q_{\ell} (\zeta), \Z/{\ell}\Z)$ corresponds to the class of
$\zeta$ under the isomorphism  $H^1 (\Q_{\ell} (\zeta),
\Z/{\ell}\Z)\cong H^1 (\Q_{\ell}(\zeta), \mu_{\ell} ) \cong
\Q_{\ell} (\zeta )^*/{\ell}$. Then the kernel of $\psi$
corresponds to $M_1$.\\

There is an unramified extension $N$ of degree ${\ell}$ over
$\Q_{\ell}$ (an Artin-Schrier extension). Let $N(\zeta ) = \Q_{\ell}
(\zeta )(\beta^{\frac{1}{{\ell}}})$ with
$\beta\in\Q_{\ell}(\zeta)^*$. Let $\varphi$ be the unramified
character in $H^1 (\Q_{\ell} , \Z/{\ell}\Z )$ whose retriction in
$H^1 (\Q_{\ell} (\zeta), \Z/{\ell}\Z)$ corresponds to the class of
$\beta$ under the isomorphism  $H^1 (\Q_{\ell} (\zeta),
\Z/{\ell}\Z)\cong H^1 (\Q_{\ell}(\zeta), \mu_{\ell} ) \cong
\Q_{\ell} (\zeta )^*/{\ell}$. Note that since $N$ is unramified,
$\beta^{\frac{1}{{\ell}}} \in  \Q_{\ell} (\zeta)^{ur}$.\\

From Tate local duality we see that $H^1 (\Q_{\ell} , \Z/{\ell}\Z )$
has the same dimension as $\Q_{\ell}^*/{\ell}$.  The latter is
isomorphic to $\Z/{\ell}\Z \oplus \Z/{\ell}\Z$.  So the dimension of
$H^1 (\Q_{\ell} , \Z/{\ell}\Z )$ is two. Since the two characters
$\psi$ and $\varphi$ are independent, one being ramified and the
other not, they form a basis of $H^1 (\Q_{\ell} , \Z/{\ell}\Z )$
over $\Z/{\ell}\Z$. It follows that every character in $H^1
(\Q_{\ell} ,\Z/{\ell}\Z)$ is of the form $a\psi + b\varphi$ with
$a,b\in \Z/{\ell}\Z$. The restriction of $a\psi + b\varphi$ in $H^1
(\Q_{\ell}(\zeta),\Z/{\ell}\Z)$ corresponds to the class of $\rho =
\zeta^a \beta^b$, and gives rise to a cyclic extension $M'$ of
degree ${\ell}$ over $\Q_{\ell}$ with $M'(\zeta) = \Q_{\ell}
(\zeta)(\rho^{\frac{1}{{\ell}}} )$. Note that $\rho \sim^{\ell}
\zeta^i$ in $\Q_{\ell} (\zeta)^{ur}$ as
$\beta^{\frac{1}{{\ell}}} \in  \Q_{\ell} (\zeta)^{ur}$.\\

Since $\varphi$ is unramified and $1+{\ell}$ is a unit,
\[ < \varphi , 1+{\ell} > = 0.\]
So
\[ < a\psi+b\varphi , 1+{\ell} > = a < \psi , 1+{\ell} > = a (\zeta , 1+{\ell} ).\]

Since $1+{\ell} = \eta_{{\ell}-1} \xi$ with $\xi\equiv 1
\pmod{\lambda^{{\ell}}}$,
$$(\eta_1, 1+{\ell})= (\eta_1 , \eta_{{\ell}-1} \xi )= (\eta_1,\eta_{{\ell}-1} )$$
$$(\eta_1 , \eta_{{\ell}-1} ) = (\eta_1 , \eta_{\ell} ) + (\eta_{\ell} , \eta_1 ) - ({\ell}-1) (\eta_{\ell} , \lambda ) = 1.$$
([CF] p.354; our symbol is written additively.)\\

Therefore, $< a\psi+b\varphi , 1+{\ell} > = a$.

The restriction of $\chi_u$ corresponds to $a\psi+b\varphi$, with
$a,b\in\Z/{\ell}\Z$, under the isomorphism between $H^1 (K_{u} ,
\Z/{\ell}\Z )$ and $H^1(\Q_{\ell},\Z/{\ell}\Z)$. From the discussion
above, $A\sim^{\ell} \zeta^a \beta^b$ under the identification of
$K(\mu_{\ell})_{u'}$ with $\Q_{\ell}(\mu_{\ell})$, and $A\sim^{\ell}
\zeta^a$ in
$\Q_{\ell}(\mu_{\ell})^{ur}$.\\

We have
\[({\ell}-1) < \chi_u  , 1+{\ell}_u  >= < \chi_{u'} , 1+{\ell}_{u'} > = < a\psi+b\varphi , 1+{\ell} > = a.\]
So
\[ n=\sigma_u(\chi)= < \chi_u , 1+{\ell}_u > = -a\]
where $A\sim^{\ell} \zeta^a$ in $\Q_{\ell}(\mu_{\ell})^{ur}$.

\vspace{.5in}

\ \\{\bf \Large{References}}
\begin{itemize}
\begin{footnotesize}
\item[{[BMSW]}]  B. Bektermirov, B. Mazur, W. Stein and M. Watkins,
\emph{Average ranks of elliptic curves:  tension between data and
conjectures}, Bull. American Math. Society 44 (2007) 233-254

\item[{[CF]}]  J.W.S. Cassels and A. Fr\"ohlich, \emph{Algebraic
Number Theory}, Academic Press 1967

\item[{[D]}]  M. Deuring, \emph{Die Typen der Multiplikatorenringe elliptischer Funktionenkörper},
Abh. Math. Sem. Hansischen Univ. 14, (1941). 197-272.

\item[{[F]}]  G. Frey, \emph{Applications of arithmetical geometry to cryptographic constructions},
 In Proceedings of the Fifth International
Conference on Finite Fields and Applications. Springer Verlag, page
128-161, 1999.
%Preprint also available
%http://www.exp-math.uni-essen.de/zahlentheorie/preprints/Index.html.

\item[{[FR]}]
G. Frey and H.-G. R\"uck, \emph{A remark concerning $m$-divisibility
and the discrete logarithm in the divisor class group of curves},
{\em Mathematics of Computation}, 62(206):865--874, 1994.

\item[{[G]}]  D. Goldfeld, \emph{Conjectures on elliptic curves over quadratic fields}, in Number Theory (Carbondale, Ill., 1979), Lecture Notes in Math. 751, Springer, Berlin, 1979, 108--118
\item[{[HB]}]  D.R. Heath-Brown, {\emph The average analytic rank of elliptic curves},
Duke Math. J. 122 (2004), no. 3, 591--623.

\item[{[HKT]}] M.-D. Huang, K. L. Kueh, and K.-S. Tan \emph{Lifting
elliptic curves and solving the elliptic curve discrete logarithm
problem} In ANTS,  Lecture Notes in Computer Science, Volume 1838
 Springer-Verlag, 2000.

%\item[{[HR1]}]  M.-D. Huang and W. Raskind, \emph{Global duality
%and the discrete logarithm problem}, preprint 2006.
%http://www-rcf.usc.edu/~mdhuang/papers.html

%\item[{[HR2]}]  M.-D. Huang and W. Raskind, \emph{Signature calculus and the
%the discrete logarithm problem for the multiplicative group case},
%preprint 2006. http://www-rcf.usc.edu/~mdhuang/papers.html

%\item[{[HR3]}]  M.-D. Huang and W. Raskind, \emph{Signature
%calculus and the discrete logarithm problem for elliptic curves},
%preprint 2006

\item[{[HRANTS]}]  M.-D. Huang and W. Raskind, \emph{Signature
calculus and discrete logarithm problems}, Proceedings of the 7th
Algorithmic Number Theory Symposium (ANTS 2006), LNCS 4076, 558-572,
Springer-Verlag, 2006.

\item[{[JKSST]}] M.J. Jacobson, N. Koblitz, J.H. Silverman, A. Stein, and E.
Teske. Analysis of the Xedni calculus attack. Design, Codes and
Cryptography, 20 41-64, 2000

\item[{[Ka]}]  S. Kamienny, \emph{Torsion points on elliptic curves and $q$-coefficients of modular forms}. Invent. Math. 109 (1992), no. 2, 221--229.

\item[{[Ko]}] N. Koblitz \emph{Elliptic curve cryptosystems} Mathematics of
Computation, 48 203-209, 1987.

\item[{[KMV]}] N. Koblitz, A. Menezes and S. Vanstone \emph{The state of elliptic
curve cryptography}, Design, Codes and Cryptography, 19, 173-193
(2000)

\item[{[Ma]}]  B. Mazur, \emph{Notes on the \'etale cohomology of number
fields}, Ann. Sci. \'Ecole Normale Sup\'erieure 6 (1973) 521-556

\item[{[Mc]}]  K. McCurley, \emph{The discrete logarithm problem},
in Cryptology and Computational Number Theory, C. Pomerance, editor,
Proceedings of Symposia in Applied Mathematics, Volume 42, 49-74,
1990

\item[{[Mill]}]  V. Miller \emph{Uses of elliptic curves in cryptography},  In
Advances in Cryptology: Proceedings of Crypto 85, Lecture Notes in
Computer Science,  volume 218, 417-426. Springer-Verlag, 1985.

\item[{[L]}]  S. Lang \emph{Algebraic groups over finite fields} Amer. J. Math. 78 (1956), 555--563.

\item[{[MET]}] J.S. Milne, \emph{\'Etale Cohomology}, Princeton Mathematical Series,
Volume 33, Princeton University Press 1980

\item[{[MAD]}]J.S. Milne, \emph{Arithmetic Duality Theorems}, Perspectives in
Mathematics, Volume 1., Academic Press 1986

\item[{[N]}]  K. Nguyen, Thesis, Universit\"at Essen, 2001

\item[{[RS]}]  K. Rubin and A. Silverberg, Torus-based cryptography,
in Advances in Cryptology --- CRYPTO 2003, Lecture Notes in Computer
Science 2729 (2003), Springer, 349-365

\item[{[S1]}]  J.-P. Serre, Corps Locaux, Paris Hermann 1962;
English translation: Local Fields, Graduate Texts in Mathematics,
Volume 67, Springer Verlag, Heidelberg-New York, 1979

\item[{[S2]}]  J.-P. Serre, \emph{Groupes p-divisibles (d'apr\`es J.
Tate)}, S\'eminaire Bourbaki 1966/67, Expos\'e 318, reprinted by the
Soci\'et\'e Math\'ematique de France 1995

\item[{[Sc]}] R. Schoof, \emph{Counting points on elliptic curves
over finite fields}, Journal de Th\'{e}orie des Nombres de Bordeaux
7 (1995), 219-254.

\item[{[Si1]}] J.H. Silverman, \emph{The
Arithmetic of Elliptic Curves}, Graduate Texts in Mathematics,
Volume 106, Springer Verlag 1986.

\item[{[Si2]}] J.H. Silverman, \emph{Advanced Topics in the
Arithmetic of Elliptic Curves}, Graduate Texts in Mathematics,
Volume 151, Springer Verlag 1994.

\item[{[SWD]}]O. Schirokauer, D.Weber, and T. Denny \emph{Discrete logarithms:
The effectiveness of the index calculus method} In ANTS II, volume
1122 of Lecture Notes in Computer Science. Springer-Verlag, 1996.

\end{footnotesize}
\end{itemize}

\end{document}